\documentclass[preprint,12pt,authoryear]{elsarticle}

\usepackage{xcolor}
\usepackage{mathrsfs}
\usepackage{lipsum}
\usepackage{amsfonts}
\usepackage{algorithmic}
\usepackage{amssymb}
\usepackage{amsmath}
\usepackage{graphicx}
\usepackage{subfig}
\usepackage{epstopdf}
\usepackage{mathtools}
\usepackage{cleveref}
\usepackage{tabularx}
\usepackage{algorithm2e}
\usepackage{wrapfig}
\usepackage{hyperref}
\usepackage{bm}
\usepackage{amsopn}

\SetKwComment{Comment}{/* }{ */}

\DeclareMathOperator{\Tr}{Tr}
\ifpdf
\DeclareGraphicsExtensions{.eps,.pdf,.png,.jpg}
\else
\DeclareGraphicsExtensions{.eps}
\fi
\newcommand{\RN}[1]{%
  \textup{\uppercase\expandafter{\romannumeral#1}}%
}
\newcommand{\rn}[1]{%
  \textup{\expandafter{\romannumeral#1}}%
}

\journal{Journal of Theoretical Biology}

\begin{document}

\begin{frontmatter}



\title{Simulating Irregular Symmetry Breaking in Gut Cross Sections Using a Novel Energy-Optimization Approach in Growth-Elasticity} 

%
\author[add1,add2]{Min Wu}
\ead{englier@gmail.com}

\address[add1]{Department of Mathematical Sciences, Worcester Polytechnic Institute}
\address[add2]{Center for Computational Biology, Flatiron Institute}
\cortext[CorrespondingAuthor]{Corresponding author}




\begin{abstract}
Growth-elasticity (also known as morphoelasticity) is a powerful model framework for understanding complex shape development in soft biological tissues. At each instant, by mapping how continuum building blocks have grown geometrically and how they respond elastically to the push-and-pull from their neighbors, the shape of the growing structure is determined from a state of mechanical equilibrium. As mechanical loads continue to be added to the system through growth, many interesting shapes, such as smooth wavy wrinkles, sharp creases, and deep folds, can form on the tissue surface from a relatively flatter geometry.

Previous numerical simulations of growth-elasticity have reproduced many interesting shapes resembling those observed in reality, such as the foldings on mammalian brains and guts. In the case of mammalian guts, it has been shown that wavy wrinkles, deep folds, and sharp creases on the interior organ surface can be simulated even under a simple assumption of isotropic uniform growth in the interior layer of the organ. Interestingly, the simulated patterns are all regular  along the tube's circumference, with either all smooth or all sharp indentations, whereas some undulation patterns in reality exhibit irregular patterns and a mixture of sharp creases and smooth indentations along the circumference. Can we simulate irregular indentation patterns without further complicating the growth patterning?

In this paper, we have discovered abundant shape solutions with irregular indentation patterns by developing a Rayleigh-Ritz finite-element method (FEM). In contrast to previous Galerkin FEMs, which solve the weak formulation of the mechanical-equilibrium equations, the new method formulates an optimization problem for the discretized energy functional, whose critical points are equivalent to solutions obtained by solving the mechanical-equilibrium equations. This new method is more robust than previous methods. Specifically, it does not require the initial guess to be near a solution to achieve convergence, and it allows control over the direction of numerical iterates across the energy landscape. This approach enables the capture of more solutions that cannot be easily reached by previous methods. In addition to the previously found regular smooth and non-smooth configurations, we have identified a new transitional irregular smooth shape, new shapes with a mixture of smooth and non-smooth surface indentations, and a variety of irregular  patterns with different numbers of creases. Our numerical results demonstrate that growth-elasticity modeling can match more shape patterns observed in reality than previously thought.
 \end{abstract}


\begin{highlights}
\item A robust energy-optimization FEM method for solving growth-elasticity models is developed.
\item The new method is implemented to simulate growth-induced gut morphogenesis.
\item Numerous irregular undulation patterns are numerically realized in gut cross sections. 
\item Irregular  patterns with mixed smooth and non-smooth indentations better match reality.
\end{highlights}


\end{frontmatter}



\section{Introduction}\label{sec:growing_tissues}
Growth plays an important role in developmental processes \citep{neufeld1998coordination,lecuit2003developmental,ingber2006mechanical}, tissue regeneration \citep{loewenstein1967intercellular,sun2014control}, cancer progression \citep{folkman2002role,altorki2019lung}, and organoid engineering \citep{lancaster2014organogenesis,hofer2021engineering}. Growth can generate forces in addition to those caused by cellular and sub-cellular contractile activities \citep{hallatschek2023proliferating}, and these growth-induced forces can be transmitted throughout the living tissues, leading to complex deformations and residual stresses \citep{nia2016solid}.

Along these lines of thought, the growth-elasticity theory, also known as morphoelasticity, has been developed \citep{Rodriguez-1994-decomposition}. For more details, see the book \citep{goriely2017mathematics}. In this theory, the incorporation of new biomass is captured geometrically by volumetric growth. Specifically, a growth-stretch tensor (matrix) field is defined to describe the change in reference size and shape of the preexisting building blocks in the growing continuum. These (soft) building blocks adjust their local configuration in response to the push-and-pull interactions with their neighbors, resulting in an elastic deformation field. With the aid of a hyperelastic strain energy function, which describes the stress-strain relationship of the tissue constituents, one can solve the global configuration(s) of the growing tissue from mechanical-equilibrium equations with proper boundary conditions.

The theory has been widely applied to understand morphogenesis in living organisms by various research groups.  For example, the optic vesicle formation in the chicken embryo \citep{garcia2017contraction,oltean2018apoptosis} and the mammal brain gyrifications \citep{budday2015size,tallinen2016growth} are explained as a result of specific patterning of growth, respectively, from numerical simulations where the initial nonlinear surface distortion exaggerates as differential and directional growth accumulates. These organ morphogenesis involve both smooth foldings and local non-smooth indentations. See also the reviews \cite{taber1995biomechanics, Ambrosi-2019-growth-review}.

In addition to simulations of morphogenesis in organs with irregular shapes, many previous studies have focused on symmetry-breaking problems in growing structures with perfect annular \citep{li2011surface,moulton2011circumferential,ciarletta2012growth,wu2015modelling} or spherical shapes \citep{goriely2005differential,amar2005growth,li2011surface2}. In these studies, symmetry breaking is a system bifurcation in which the annular or spherical solutions become unstable, giving rise to new solutions with more complex geometries. The symmetry breaking of the annular shape at the interior surface has significant biological implications for gut crypt formation (a process initialized by indentations and bulges on the epithelial layer)  and tumorigenesis, leading to numerous studies on this geometric setting \citep{li2011surface,jin2011creases,moulton2011circumferential,ciarletta2012growth,papastavrou2013mechanics,balbi2015morphoelastic}.

Since the gut tubes are layered with different cells and extracellular-matrix constituents that indicate different growth activities and different material properties among the layers, previous works reveal these two ingredients are sufficient to generate folding patterns that resemble the gut interior surface morphologies. In a simplified bilayer setting, they have shown that the interior-boundary folding patterns can emerge due to a combination of growth concentrated in the interior layer and the geometric confinement from the exterior boundary surface. Simulations show that smooth wave patterns can emerge when the growing layer is much stiffer than the non-growing layer \citep{li2011surface,papastavrou2013mechanics}; and some also show that smooth indentations can transition into sharp creases all at once \citep{balbi2015morphoelastic}. In the extreme setting where the external layer is rigid, patterns with regular sharp creases along the interior surface are formed instead of smooth wrinkles or foldings \citep{jin2011creases}.

\begin{wrapfigure}{l}{0.5\textwidth}
\includegraphics[width=0.9\linewidth]{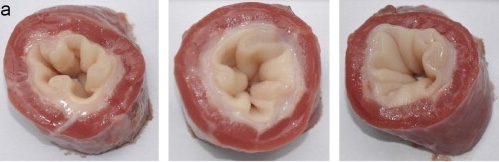} 
\caption{The bovine esophagus cross sections,  reprinted from \cite{li2011surface} with permission from Elsevier. The interior surfaces exhibit irregular patterns and a mixture of sharp creases and smooth indentations along the circumference.}
\label{Fig-1}
\end{wrapfigure}

Interestingly, these simulated morphologies are all regular along the tube's circumference, featuring either all smooth or all sharp indentations. However, some undulation patterns in reality exhibit irregular patterns and a mixture of sharp creases and smooth indentations along the circumference. For example, see Fig.\ref{Fig-1}. Can we simulate irregular  patterns without complicating the growth pattern assumptions within the framework of growth-elasticity?

In this paper, we have discovered multiple irregular  patterns in a bi-layer annular geometry by developing a Rayleigh-Ritz finite-element method (FEM). Most previous FEM methods are Galerkin FEMs, which solve the mechanical-equilibrium equations from growth-elasticity without fully exploiting the elastic energy landscape \citep{li2011surface,jin2011creases,papastavrou2013mechanics,balbi2015morphoelastic}. These methods are typically constrained by the regular perturbations of the symmetric annular solution, followed by Newton's iteration to approach new shape solutions, provided they are sufficiently close to the initial guess after the perturbation. In contrast, our Rayleigh-Ritz method is based on the optimization of the discretized energy functional, with its critical points corresponding to solutions derived from solving the mechanical-equilibrium equations.

This new solution-search strategy is more robust than previous methods. Specifically, it does not require the initial guess to be near a solution to achieve convergence, and it allows control over the direction of numerical iterates across the energy landscape. This approach enables the capture of more solutions that cannot be easily reached by Newton's iterations alone, without needing to implement pseudo-arclength numerical continuation \citep{groh2022morphoelastic}. It successfully finds irregular  shape solutions that qualitatively match experimental observations \citep{li2011surface}, which have not been numerically realized before. While the Rayleigh-Ritz method is well-known for solving nonlinear elastic problems \citep{reddy2017energy}, there have been few developments in applying it to growth-elasticity models \citep{tallinen2013surface,tallinen2016growth}. As an improvement to \cite{tallinen2013surface,tallinen2016growth}, we also leverage Hessian information to design stable gradient-descent iterations and perform local stability analysis of the numerical solutions.

\begin{wrapfigure}{l}{0.5\textwidth}
\includegraphics[width=0.9\linewidth]{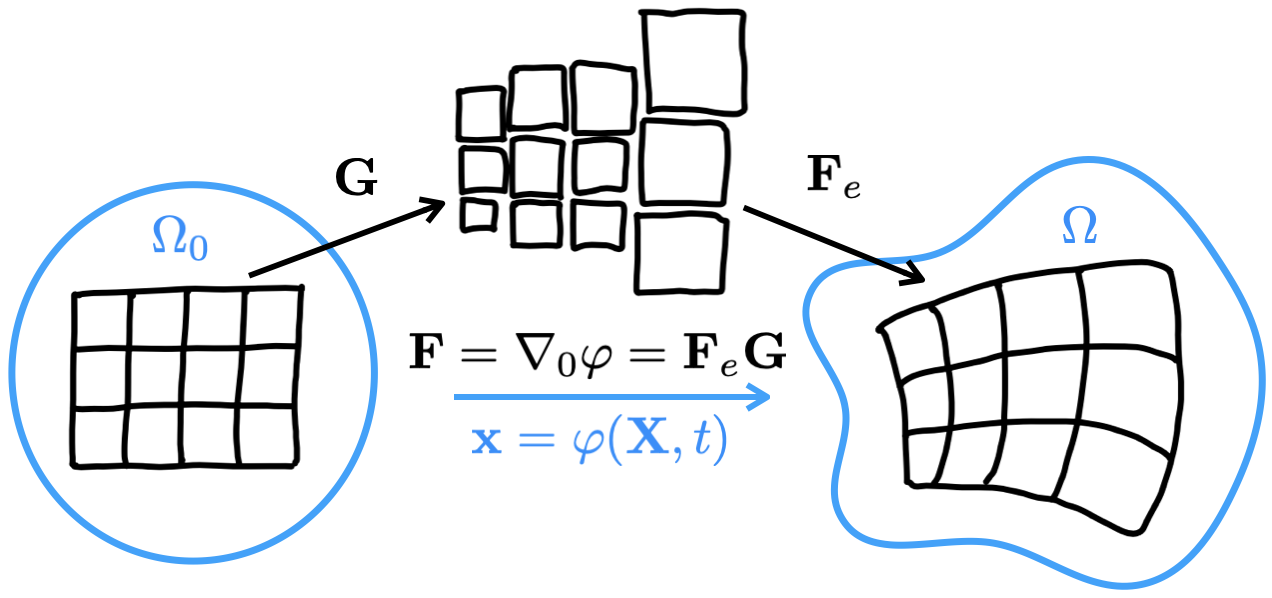} 
\caption{The transformation of tissue from configuration $\Omega_0$ to $\Omega$ is captured by the mapping $\mathbf{x}=\varphi(\mathbf{X},t)$, where its spatial gradient $\nabla_0\varphi:=\mathbf{F}$ defines  local geometric deformation of tissue particles/blocks. The deformations are explained by the distribution of growth $\mathbf{G}$ and elastic deformations $\mathbf{F}_e$ through the decomposition $\mathbf{F}=\mathbf{F}_e \mathbf{G}$. }
\label{Fig0}
\end{wrapfigure}


The paper is organized as follows. In Sec. 2, we present the growth-elasticity model and the 2D bilayer annular growth problem we aim to solve, which simulates the shape development of the gut-tube cross sections.  In Sec. 3, we present the formulation of the new Rayleigh-Ritz FEM method and the algorithm. In Sec. 4, we apply the algorithm to study three different cases of the bilayer problem, where the material stiffness ratio between the growing and non-growing layers is varied. Among the three cases,  we identify a variety of different shapes. Besides the regular smooth and  non-smooth configurations that have been found before, we find a new transitory irregular  smooth shape, new shapes with a mixture of smooth and non-smooth surface indentations, and a bountiful irregular non-smooth  patterns with different number of creases.  In Sec. 5, we discuss our results with previous works and shape patterns of gut tubes in reality, and suggest further modifications of growth-elasticity theory given the insight from our numerical solutions.

\section{The model formulation}\label{sec:computation}
In this section, we present growth-elasticity theory in its variational formulation in conjunction with its Euler-Lagrangian equations that gives the nonlinear partial differential equations system. We will begin by introducing the multiplicative decomposition of the deformation gradient, and end with the specific mathematical problem we aim to solve in 2D. Our presentation is in the Lagrangian coordinate system.

\subsection{Decomposition of deformation gradient}
In growth-elasticity (see Fig.\ref{Fig0}), one follows the Kr\"oner-Lee approach \cite{kroner1959allgemeine,lee1969elastic} to decompose the deformation gradient $\mathbf{F}:= {\partial \mathbf{x}}/{\partial \mathbf{X}}=\nabla_0\varphi$ by
\begin{eqnarray}\label{eq_decomp}
\mathbf{F}= \mathbf{F}_e \mathbf{G}
\end{eqnarray}
where $\varphi(\mathbf{X},t):=\mathbf{x}$ is the observable motion, mapping a particle at $\mathbf{X}$ initially to its current position $\mathbf{x}$. $\mathbf{F}_e(\mathbf{X},t)$ is the elastic deformation tensor associated with stress. Different from the Kr\"oner-Lee approach where $\mathbf{G}$ is considered as a plastic deformation, here, $\mathbf{G}(\mathbf{X},t)$ is the accumulated growth (stretch) tensor such that $J_g= \det{\mathbf{G}}>0$ describes the local growth in volume. In particular, $J_g>1$ corresponds to a net growth and $J_g<1$ corresponds to a net atrophy. Since $\det{\mathbf{G}}>0$, one can apply polar decomposition on  $\mathbf{G}=\mathbf{U}\mathbf{Q}=\mathbf{Q}'\mathbf{V}$, where the volumetric growth is interpreted by a positive-definite growth stretch tensor $\mathbf{U}$ ($\mathbf{V}$) specifying the change of length along orthogonal directions preceded (followed) by a local rotation $\mathbf{Q}$ ($\mathbf{Q}'$).
\subsection{The variational formulation of growth-elasticity theory}
When the trajectory of $\mathbf{G}(\mathbf{X},t)$ is prescribed, one can solve  $\varphi(\mathbf{X},t)$ and $\mathbf{F}_e=\mathbf{F}\mathbf{G}^{-1}=\nabla_0\varphi \mathbf{G}^{-1}$ through the optimization of the energy functional \cite{ciarletta2012growth,Wu-2015-wound}
\begin{eqnarray}
E=\int_{\Omega_0} W(\mathbf{F}_e) J_g d \mathbf{X}-\int_{\partial\Omega_0^N} \mathbf{T}\cdot \varphi d\mathbf{S}
\label{energy_old}
\end{eqnarray} 
where $\Omega_0\subset R^d$ ($d=2,3$) is the tissue domain in the Lagrangian frame, $W=W(\mathbf{F}_e)$ is a strain-energy density function, and $J_g$ works as the ``Jacobian" from initial Lagrangian coordinate $\mathbf{X}$ to a stress-free grown state. The boundary integral term specifies the energy potential due to some traction field $\mathbf{T}$ along the boundary $\partial\Omega_0^N$. We note that $\mathbf{T}$ here does not depend on $\varphi$ or its derivatives -- otherwise it is a non-trivial task to find the boundary potential \cite{podio1991surface}. We also note that when $\mathbf{G}=\mathbf{I}$ such that $J_g=1$, the growth-elasticity system recovers the purely hyperelastic framework with $\mathbf{F}_e=\mathbf{F}$.

The Euler–Lagrange equation of Eq.(\ref{energy_old}) gives the standard mechanical equilibrium 
 \begin{eqnarray}
\begin{cases}
    \nabla_0\cdot \boldsymbol{\Pi} &= \mathbf{0},  \hspace{1cm}\boldsymbol{\Pi} = J_g\frac{\partial W}{\partial \mathbf{F}_e}\mathbf{G}^{-T} \text{ in }\Omega_0\\
 \boldsymbol{\Pi}\mathbf{N} &= \mathbf{T} \text{ on }\partial\Omega_0^N
\end{cases}
\label{stationary}
\end{eqnarray}
with  $\boldsymbol{\Pi}$ the first Piola–Kirchhoff stress tensor, related to the Cauchy stress by the Piola transform $\boldsymbol \sigma = (\det\mathbf{F})^{-1}\boldsymbol{\Pi}\mathbf{F}^T$, and  $\nabla_0$ the gradient operator in the Lagrangian coordinate, and $\mathbf{N}$ is the outward normal on moving boundary $\partial\Omega_0^N$. The two sub-equations correspond to the force balance in the bulk and along the moving boundary, respectively. 


\subsection{The bilayer annular growth problem}
\label{sec:model}

\begin{wrapfigure}{l}{0.3\textwidth}
\includegraphics[width=1.0\linewidth]{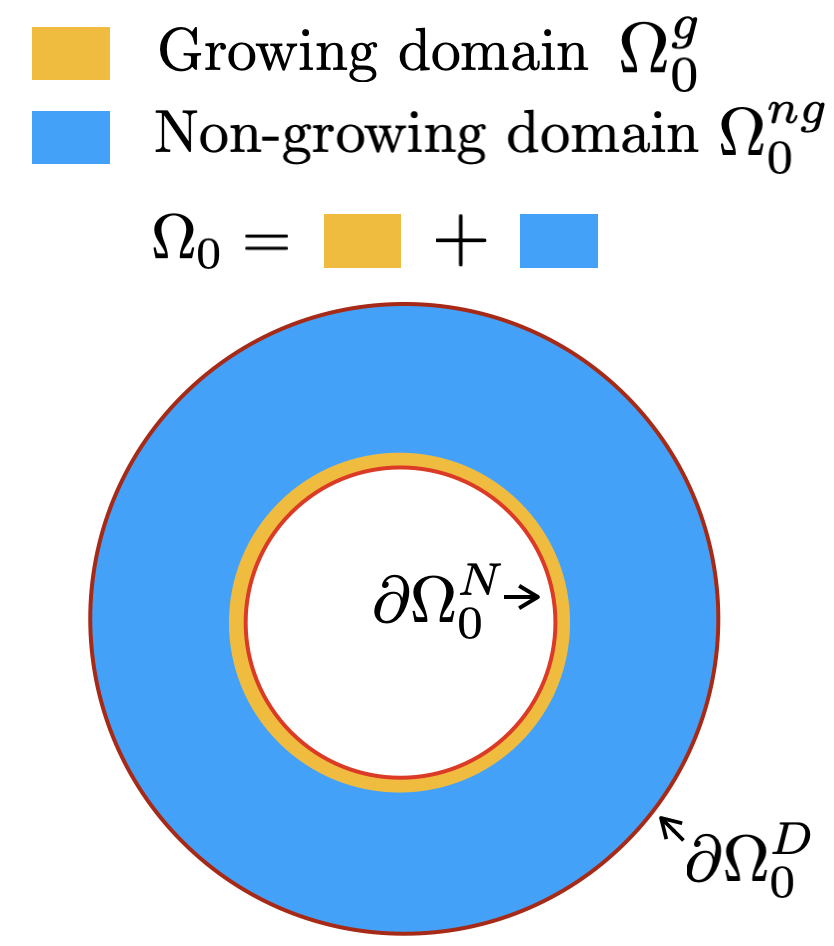} 
\caption{Schematics of the initial domain.}
\label{Fig2}
\end{wrapfigure}

In this paper, we simulate the growth-driven shape dynamics of a bilayer annulus, which resembles the 2D cross section of gut tubes (see Fig.\ref{Fig2}). Previous works concerning this problem solve Eq.(\ref{stationary}) with the multiplicative decomposition Eq.(\ref{eq_decomp}) through its weak formulation \citep{li2011surface,jin2011creases,papastavrou2013mechanics,balbi2015morphoelastic}. Instead, we will develop an energy-optimization approach that solves Eq.(\ref{stationary})  via finding the critical points of the energy Eq.(\ref{energy_old}). We are particularly interested in finding new solutions, thus new phenomena,  that have not been revealed by previous approaches.  
For simplicity, we consider the interior boundary $\partial \Omega_0^N$ is free and $\mathbf{T}=\mathbf{0}$ for $\mathbf{X}\in\partial \Omega_0^N$, which reduces Eqs. (\ref{eq_decomp}) and (\ref{energy_old}) to
\begin{eqnarray}
E=\int_{\Omega_0} W(\mathbf{F}_e) J_g d \mathbf{X}=\int_{\Omega_0} W(\nabla_0\varphi \mathbf{G}^{-1}) J_g d \mathbf{X} \label{energy_simplified}
\end{eqnarray} 
in its energy functional. The second equality used Eq. (\ref{eq_decomp}). We consider the trajectory of $\mathbf{G}(\mathbf{X},t)$ as given, and the solution field $\varphi(\mathbf{X},t)$ is a minimizer of the energy Eq.(\ref{energy_simplified}) with the material points on the exterior boundary $\partial\Omega_0^D$ being fixed:
\begin{eqnarray}\varphi(\mathbf{X},t) = \mathbf{X}, \text{ for }\mathbf{X}\in\partial \Omega_0^D. \end{eqnarray} Thus, the positions $\varphi(\mathbf{X},t)$ in the bulk $\Omega_0$ and on free boundary $\Omega_0^N$ will be solved through the energy minimization.

Previous works have shown that by only having differential isotropic growth between the two layers is sufficient to induce symmetry breaking of the annular configuration \citep{li2011surface,jin2011creases,papastavrou2013mechanics}. Thus, we consider   \begin{eqnarray}
\mathbf{G}(\mathbf{X},t)=\begin{cases}
    g(t)\mathbf{I} & \text{ for } \mathbf{X}\in\Omega_0^g\\
 \mathbf{I} & \text{ for } \mathbf{X}\in\Omega_0^{ng}
\end{cases}
\label{growth_rule}
\end{eqnarray}
 where the growth only happens in the interior layer (see the yellow region in Fig.\ref{Fig2}) for simplicity. The identity matrix $\mathbf{I}$ ensures the isotropy of the growth. In this work, since the focus is on the 2D simulation, we will consider $\mathbf{I}$ as a $2\times2$ identity matrix, denoted by $\mathbf{I}_2$. Although the exterior layer does not grow, it can deform in response to the growth of the interior layer. 

When the growth function  $g(t)$ is specified, as long as it increases monotonically(such as in \citep{papastavrou2013mechanics}), it is a bijection relating the time with the accumulated growth stretch $g$. As such, we ignore the speculation on the form of $g(t)$ and treat $g$ as the single parameter of the growth-elasticity model. In particular, we study the evolution of shape solutions from the model as $g$ increases from unity in the growing domain $\Omega_0^g$. 

Previous works have shown that the different material properties between the growing and the non-growing layers can result in different shape patterns \citep{li2011surface,jin2011creases,papastavrou2013mechanics,balbi2015morphoelastic}.  For simplicity, we choose a neo-Hookean strain energy  \cite{holzapfel2002nonlinear}
\begin{eqnarray} 
W(\mathbf{F}_e)={\mu}/{2}\big((\det\mathbf{F}_e)^{-1}\Tr(\mathbf{F}_e^T\mathbf{F}_e)-2\big)+{K}/{2}\big(\det\mathbf{F}_e-1\big)^2
\label{neo}
\end{eqnarray} 
and adjust the ratios in $K$ and $\mu$, the bulk modulus and the shear modulus, respectively, to compare with previous results.  In particular, we define $K_{g}$ ($\mu_g$) and $K_{ng}$ ($\mu_{ng}$) to be the bulk (shear) modulus in the growing and non-growing region, respectively. In the results, we adjust the ratios $\mu_{ng}/\mu_{g}$ and $K_{ng}/K_{g}$ simultaneously to represent the change of stiffness ratios between the two layers.

\section{The numerical algorithm}\label{sec:computation}
\begin{wrapfigure}{l}{0.25\textwidth}
\includegraphics[width=0.9\linewidth]{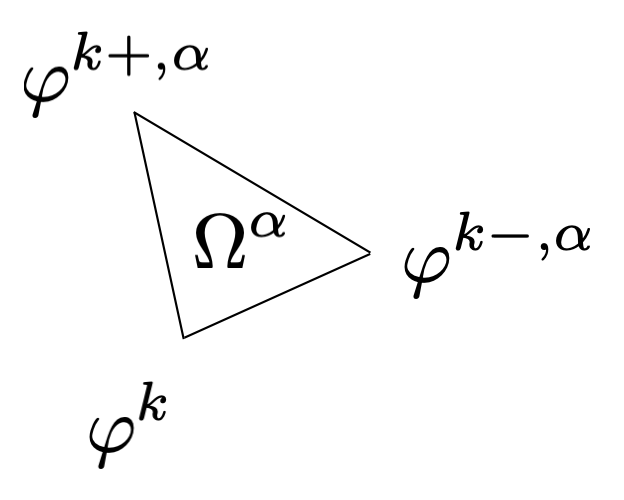} 
\caption{$\varphi^{k}\to \varphi^{k-,\alpha_k}\to\varphi^{k+,\alpha_k}$ are in the counterclockwise order in each triangle $\alpha$.}
\label{Fig1}
\end{wrapfigure}

In this section, we present our Rayleigh-Ritz FEM to solve the growth-elasticity model presented in Sec.\ref{sec:model}. Although we will only study isotropic $\mathbf{G}(\mathbf{X},t)=g\mathbf{I}_2$ treating $g$ as a control parameter, we will present the method with the form of $\mathbf{G}(\mathbf{X},t)$ unspecified. At each time point $t_j\in[0,T_{end}]$, we solve $\varphi(\mathbf{X},t_j)$ from minimizing Eq.(\ref{energy_simplified}), with $\mathbf{G}(\mathbf{X},t_j)$ given. We provide details on the spatial discretization of $\varphi(\mathbf{X},t_j)$ and its associated energy function, followed by the gradient and Hessian derivations. In the end, we provide details on the optimization algorithm, and how we apply it to solve the system in time with bifurcation analysis and numerical continuation.

\subsection{Uniform-strain Energy Discretization}
We consider the initial domain $\Omega_0=\Omega(t_0)=\cup_{\alpha}\Omega^\alpha(t_0)\in R^2$ being triangulated with nodes $\{\mathbf{X}^k\}$, where $k$ indexes all the nodes. Given the triangulation, we consider the approximation of $\varphi$ using linear continuous interpolation functions such that $\varphi$ is determined by $\{\varphi^k(t_j)\}$, and $\varphi^k(t_0)=\mathbf{X}^k$.
Thus, the deformation gradient $\mathbf{F}^\alpha=[\nabla_0{\varphi}]^\alpha$ is uniform in each $\Omega^\alpha$, determined by $\{\varphi^k(t_j)\}$ and $\varphi^k(t_0)$.  In particular, on each triangle $\Omega^\alpha$ (see Fig. \ref{Fig1}), the deformation gradient is computed from $\mathbf{F}^\alpha=\mathbf{m}^\alpha [\mathbf{M}^{\alpha}]^{-1}$ where the two matrices $\mathbf{m}^\alpha=[\varphi^{k-,\alpha}(t_j)-\varphi^{k}(t_j),  \varphi^{k+,\alpha}(t_j)-\varphi^{k}(t_j)]$ and $\mathbf{M}^\alpha=[\varphi^{k-,\alpha}(t_0)-\varphi^{k}(t_0),  \varphi^{k+,\alpha}(t_0)-\varphi^{k}(t_0)]$. 

Within each triangle $\Omega^\alpha$, we define the average elastic deformation $\mathbf{F}_e^\alpha = \mathbf{F}^\alpha [\mathbf{G}^{\alpha}]^{-1}$ given $\mathbf{G}^{\alpha}$ being the averaged $\mathbf{G}(\mathbf{X},t_j)$ in $\Omega^\alpha$.
Given the triangulation and the uniform-strain setting, we discretize the energy functional Eq.(\ref{energy_simplified}) by 
\begin{eqnarray}\label{Energy_d} E&=&
\sum_\alpha E^\alpha \text{ with }  E^\alpha(\varphi^{1,\alpha},\varphi^{2,\alpha},\varphi^{3,\alpha}, \mathbf{G}^\alpha) := J_g^\alpha W^\alpha A_0^\alpha\end{eqnarray}
where $J_g^\alpha =\det \mathbf{G}^\alpha$, $W^\alpha = W(\mathbf{F}_e^\alpha)=W(\mathbf{m}^\alpha [\mathbf{M}^{\alpha}]^{-1}[\mathbf{G}^\alpha]^{-1})$, and $A_0^\alpha=\frac{1}{2}\det\mathbf{M}^\alpha$. 

For nodes on the fixed boundary $\varphi^{k}(t_0)\in \partial\Omega_0^D$, we have $\varphi^{k}(t_j) = \varphi^{k}(t_0)$. To summarize, with $\mathbf{G}^{\alpha}$ given, we write $E(\varphi^1,\varphi^2,...,\varphi^N)$ as a multivariable function of the coordinates of all the free nodes (those in the bulk and along the free boundary). From now on, we define the global vector $\boldsymbol{\varphi} = [\varphi^1,\varphi^2,...,\varphi^N]$ concatenating the $x-$ and $y-$ components of all the free $\varphi^k$'s.

 \subsection{Gradient and Hessian components}
A critical point of the energy $E(\boldsymbol{\varphi})$ must satisfy $\partial E/\partial \boldsymbol{\varphi} =\mathbf{0}$  where $\partial E/\partial \boldsymbol{\varphi}$ is the gradient vector of the energy with respect to the global vector. Here we provide the analytical formula of $\frac{\partial E}{\partial \varphi^k}\in R^2$ on each node $k$ (see more details of the derivations in the Appendix A):
\begin{eqnarray}
\frac{\partial E}{\partial \varphi^k}&=&\sum_{\alpha_{k}}\frac{\partial E^{\alpha_k}}{\partial \varphi^k}
\end{eqnarray}
where $\{\alpha_k\}$ comprises all triangles connected to node $k$ and each
\begin{eqnarray}
\frac{\partial E^{\alpha_k}}{\partial \varphi^k}&=&(J_g^{\alpha_k}  \underbrace{\frac{\partial W^{\alpha_k}}{\partial \mathbf{F}_e^{\alpha_k}} [\mathbf{F}_e^{\alpha_{k}}]^{T}}_{J_e^{\alpha_k}\boldsymbol{\sigma}^{\alpha_k}}:\frac{\partial \mathbf{m}^{\alpha_k}}{\partial \varphi^k} [\mathbf{m}^{\alpha_k}]^{-1})A_0^{\alpha_k}
  =\frac12\boldsymbol{\sigma}^{\alpha_k}\varphi^{k\perp,\alpha_k}\label{force}
\end{eqnarray}
with $\boldsymbol{\sigma}^{\alpha_k}= (J_e^{\alpha_k})^{-1}[{\partial W^{\alpha_k}}/{\partial \mathbf{F}_e^{\alpha_k}} ][\mathbf{F}_e^{\alpha_{k}}]^{T}$ the Cauchy stress in the triangle $\alpha_k$.
The notation $\varphi^{k\perp,\alpha_k}$ is the counterclockwise orthogonal rotation of the edge vector $(\varphi^{k+,\alpha_k}-\varphi^{k-,\alpha_k})$ opposing node $k$ in the triangle $\alpha_k$. We note this equation is equivalent to the two residual equations on each node from the Bubnov–Galerkin method using linear continuous basis/interpolation functions with triangle meshes. Interestingly, this formula is more straightforward to interpret the local discretized mechanics: the force $\vec{F}_{\alpha_k}$ on each node $k$ from one connected triangle $\alpha_k$ is $\vec{F}_{\alpha_k}=-\frac{\partial E^{\alpha_k}}{\partial \varphi_k}=-\frac12\boldsymbol{\sigma}^{\alpha_k}\varphi^{k\perp,\alpha_k}=\frac{1}{2}|\varphi^{k\perp,\alpha_k}|[\boldsymbol{\sigma}^{\alpha_k}\mathbf{n}]$ where $\mathbf{n} = -\varphi^{k\perp,\alpha_k}/|\varphi^{k\perp,\alpha_k}|$ is the outward unit normal along the apposing edge of node $k$, and $\frac{1}{2}|\varphi^{k\perp,\alpha_k}|$ is half of the edge length. In mechanical equilibrium, each node $k$ satisfies $\sum_{\alpha_{k}}\vec{F}_{\alpha_k}=\vec{0}$.

We also need the Hessian matrix $H:=\partial^2 E/(\partial \boldsymbol{\varphi}\partial \boldsymbol{\varphi})$ to check if a critical point is a local minimizer or a saddle point. The Hessian matrix is equivalent to the stiffness matrix in previous Galerkin-FEMs, which seemed to be computed by symbolic operations or numerical differentiation of the residual equations, since they have not been provided explicitly. It is sparse such that the non-zero entries only come from the second-order derivatives of $E$ with the same node, or the mixed second-order derivatives involving two nodes on the same triangle.  Here we provide the analytics of Hessian components involving two nodes $k_1$ and $k_2$ (including the case when $k_1=k_2$) via triangle $\alpha$, which we use to assemble $H$. Since both $x-$ and $y-$components of $\varphi^{k_1}$ and $\varphi^{k_2}$ are involved, we can write them in terms of a $2\times2$ submatrix. For simplicity, we only compute the submatrix corresponding to the neo-Hookean energy Eq.(\ref{neo}):
 \begin{eqnarray}
 \label{total_H}
\bigg[\frac{\partial^2 E^\alpha}{\partial \varphi^{k_1}\partial \varphi^{k_2}}\bigg]_{2\times2}=\frac{\mu}{2}\bigg[\frac{\partial^2 E_\mu^\alpha}{\partial \varphi^{k_1}\partial \varphi^{k_2}}\bigg]_{2\times2} +\frac{K}{2}\bigg[\frac{\partial^2 E^\alpha_K}{\partial \varphi^{k_1}\partial \varphi^{k_2}}\bigg]_{2\times2}\label{Hessian}
\end{eqnarray}
where 
 \begin{eqnarray}
\bigg[\frac{\partial^2 E^\alpha_\mu}{\partial \varphi^{k_1}\partial \varphi^{k_2}} \bigg]=\frac{(\varphi^{k_1\perp}\cdot \varphi^{k_2\perp})}{\det[\mathbf{m}^\alpha]}[ \mathbf{B}^\alpha]^{-1}-\frac{\text{tr}[\mathbf{B}^\alpha]\epsilon(\varphi^{k_1},\varphi^{k_2})}{2\det{[\mathbf{B}^\alpha]}}\begin{bmatrix}
0 & -1\\
1 & 0 
\end{bmatrix}
\end{eqnarray}
and 
 \begin{eqnarray}
\bigg[\frac{\partial^2 E^\alpha_K}{\partial \varphi^{k_1}\partial \varphi^{k_2}}\bigg] =\frac{1}{J_g^\alpha\det[\mathbf{M}^\alpha]}[\varphi^{k_1\perp}\otimes \varphi^{k_2\perp}]+(J_e^\alpha -1)\epsilon(\varphi^{k_1},\varphi^{k_2})\begin{bmatrix}
0 & -1\\
1 & 0 
\end{bmatrix}
\end{eqnarray}

with $[\mathbf{B}^\alpha] = [\mathbf{F}_e^\alpha][\mathbf{F}_e^\alpha]^T$ and $\epsilon(\varphi^{k_1},\varphi^{k_2})=\left\{ 
  \begin{array}{ c l }
   0 & \quad \textrm{if } \varphi^{k_1} =\varphi^{k_2} \\
    -1                 & \quad \textrm{if }\varphi^{k_1} =\varphi^{k_2+,\alpha}\\ 
     1                 &\quad \textrm{if } \varphi^{k_1} =\varphi^{k_2-,\alpha}
  \end{array}.\right.$

 \subsection{Iterative Optimizations}
 \label{algo}
We solve $\boldsymbol{\varphi}(t)$ from optimizing the energy $E(\boldsymbol{\varphi},t)$ for a set of discrete time points $\{t_0,t_1,...,t_N\}$. The energy function $E(\boldsymbol{\varphi},t)$ is not convex in general, as multiple critical points are expected. As such, the solution $\boldsymbol{\varphi}(t_j)$  we approach at time $t_j$ depends on the initial guess $\boldsymbol{\varphi}_0$. We will describe at the end of the subsection how we choose $\boldsymbol{\varphi}_0$ differently to fulfill different solution-search purposes. To begin, we present the iterative algorithm when  $\boldsymbol{\varphi}_0$ is given. There are two stages of the iterations, first by incremental energy-stable steepest gradient-descent iterations, then Newton's iterations.  
As the standard, the steepest gradient-descent iterates follow
  \begin{eqnarray}
\boldsymbol{\varphi}_{i+1}-\boldsymbol{\varphi}_{i} = -\Delta s \partial E/\partial \boldsymbol{\varphi}_i
  \end{eqnarray}
  and the 
  Newton's iterates follow
  \begin{eqnarray}
\boldsymbol{\varphi}_{i+1}-\boldsymbol{\varphi}_{i} = -[\partial^2 E/(\partial \boldsymbol{\varphi}_i\partial \boldsymbol{\varphi}_i)]^{-1}\partial E/\partial \boldsymbol{\varphi}_i.
  \end{eqnarray}
 For the gradient-descent iterations, we only take incremental steps $\Delta s$ such that the displacement between successive iterates  $||\boldsymbol{\varphi}_{i+1}-\boldsymbol{\varphi}_{i}||_2\ll1$ is small.  The assumption of incremental displacement allows us to derive the following step-size conditions to guide the iterations. The first conditions is to facilitate energy stability $E(\boldsymbol{\varphi}_{i})\geq E(\boldsymbol{\varphi}_{i+1})$:
 \begin{eqnarray}
\Delta s\leq\Delta s_{stable}:={2}/{\lambda_{max}[\partial^2 E/(\partial \boldsymbol{\varphi}_i\partial \boldsymbol{\varphi}_i)]} \text{ when } {\lambda_{max}[\partial^2 E/(\partial \boldsymbol{\varphi}_i\partial \boldsymbol{\varphi}_i)]}>0
\label{stable}
\end{eqnarray}
where  $\lambda_{max}[\Box]$ denotes the largest eigenvalue of $[\Box]$. One can see a potential problem of choosing a meaningful $\Delta s_{stable}$ if ${\lambda_{max}[\partial^2 E/(\partial \boldsymbol{\varphi}_i\partial \boldsymbol{\varphi}_i)]}\leq 0$, which means the Hessian is negative semi-definite.  In this case, any incremental step size $\Delta s$ that ensure an incremental displacement $||\boldsymbol{\varphi}_{i+1}-\boldsymbol{\varphi}_{i}||_2\ll1$ will guarantee the energy stability. See more details in the Appendix B. In fact, we have not encountered this extreme case and have only seen ${\lambda_{max}[\partial^2 E/(\partial \boldsymbol{\varphi}_i\partial \boldsymbol{\varphi}_i)]}>0$ in practice. 

When the gradient-descent iterates are close to a solution with Hessian being positive semi-definite (${\lambda_{min}[\partial^2 E/(\partial \boldsymbol{\varphi}_i\partial \boldsymbol{\varphi}_i)]}\geq 0$), we also derive the optimal step size which gives the fastest rate of 1st-order convergence: \begin{eqnarray}\Delta s=\Delta {s}_{opt} := {2}/\big(\lambda_{max}[\partial^2 E/(\partial \boldsymbol{\varphi}_i\partial \boldsymbol{\varphi}_i)]+ \lambda_{min}[\partial^2 E/(\partial \boldsymbol{\varphi}_i\partial \boldsymbol{\varphi}_i)]\big)\label{optimal}\end{eqnarray} 
See derivation in the Appendix B. In practice, we choose $\Delta {s} = \min{(\Delta {s}_{opt}, \Delta s_{stable})}$. The solution-search protocol is summarized in Algorithm 1.  \begin{algorithm}
\label{algo}
\caption{An iterative algorithm at each time $t_j$.  See texts for the choice of $\boldsymbol\varphi_0$, $Tol$, $MaxIter_1$, and $MaxIter_2$. }\label{alg:two}
\KwData{$[\mathbf{G}]^\alpha(t_j)$ and $\boldsymbol\varphi_0$. }
\KwResult{$\boldsymbol\varphi(t_j)$.}
Start with the initial guess $\boldsymbol\varphi_0$;

\While{$||\partial E/\partial \boldsymbol{\varphi}_i ||_2 > Tol$ and $Iter<MaxIter_1$}{
$\boldsymbol{\varphi}_{i+1}=\boldsymbol{\varphi}_{i}  -\Delta s \partial E/\partial \boldsymbol{\varphi}_i$ with $\Delta {s} = \min{(\Delta {s}_{opt}, \Delta s_{stable})}$;
}

\While{$||\partial E/\partial \boldsymbol{\varphi}_i ||_2 > Tol$ and $Iter<MaxIter_2$}{
$\boldsymbol{\varphi}_{i+1}=\boldsymbol{\varphi}_{i} -[\partial^2 E/(\partial \boldsymbol{\varphi}_i\partial \boldsymbol{\varphi}_i)]^{-1}\partial E/\partial \boldsymbol{\varphi}_i.$
}
\end{algorithm}

We will apply Algorithm 1 in different ways to find shape solutions. We summarize them into the three following cases: 1) In the numerical-continuation study of a solution dynamics $\boldsymbol\varphi(t)$ as time increases, the initial guess at each time point $t_j$ is given by $\boldsymbol\varphi_0 =\boldsymbol\varphi(t_{j-1})$, the solution from the last time step. 2) In the case when $\boldsymbol\varphi(t_j)$ is a saddle point at time $t_j$, we may be interested in finding a new critical point $\tilde{\boldsymbol\varphi}(t_j)$, especially a new local minimizer. In this case, we use the initial guess  $\boldsymbol\varphi_0=\boldsymbol\varphi(t_j)+\vec{v}$ where $\vec{v}$ is an eigenvector of the Hessian at $\boldsymbol\varphi(t_j)$. In most cases, the eigenvector $\vec{v}$ has a negative eigenvalue, such that the initial guess nudges $\boldsymbol\varphi(t_j)$ to a new position with a lower energy, based on which the gradient-descent iterates can approach a new critical point with a lower energy. 3) In the case when a new solution $\tilde{\boldsymbol\varphi}(t_j)$ is found, we are interested in checking if this new solution branches from the old solution $\boldsymbol\varphi(t_j)$ or was born subcritically at an earlier time. Thus, we perform a numerical-continuation study backward in time until a threshold when $\tilde{\boldsymbol\varphi}(t)$ disappears. In this case, the initial guess at each time point $t_{j'}$ is given by $\tilde{\boldsymbol\varphi}_0 =\tilde{\boldsymbol\varphi}(t_{j'+1})$.

We set error tolerance to be $Tol=10^{-7}$.  In most results, the numerical solutions meet the error tolerance under $MaxIter_1 = 10^4$ and $MaxIter_2 = 30$. For particularly challenging cases, we need to adjust $MaxIter_1=10^5-10^6$ (e.g., Fig.\ref{Fig_revision}).  

\section{Results}
\begin{figure}[!b]
	\centering
	\includegraphics[width=1.0\textwidth]{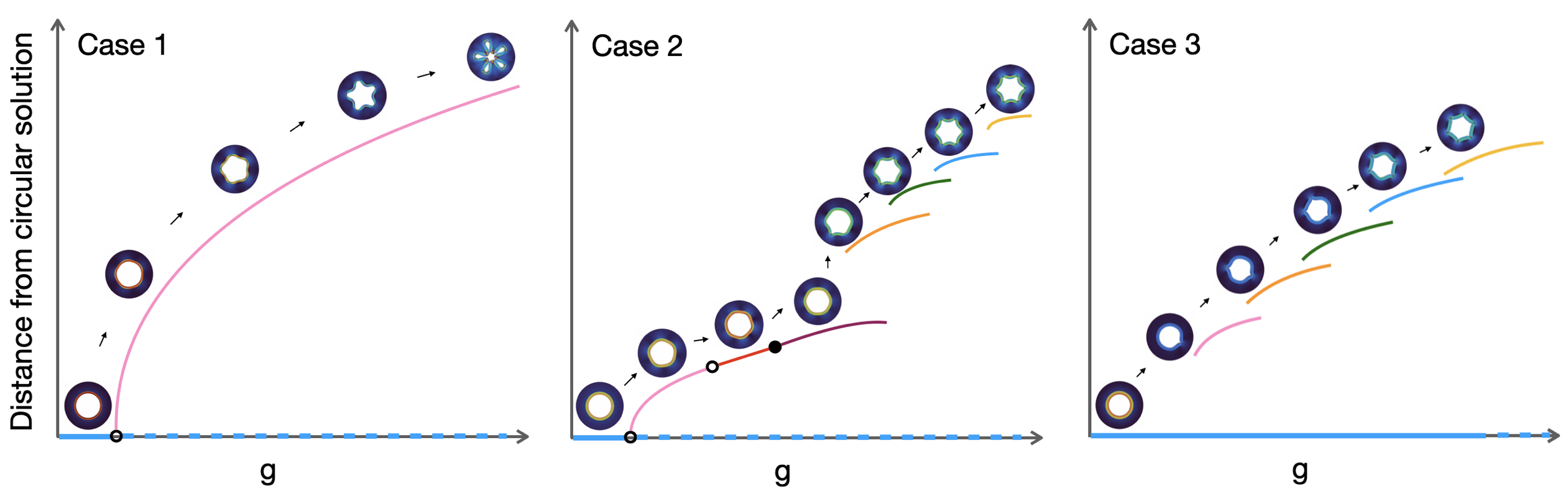}
	\caption{Summary diagrams of the morphological development in the three study cases. In case 1 with $\mu_{ng}/\mu_{g}=0.1$ and $K_{ng}/K_{g}=0.1$ (left), regular smooth wavy wrinkles branch (see the empty circle) from the circular configuration when it becomes locally unstable, which further develop into deep folds. In case 2 with $\mu_{ng}/\mu_{g}=0.4$ (middle) and $K_{ng}/K_{g}=0.4$, regular smooth wavy wrinkles branch from the circular configuration initially (see the first empty circle) but cannot further develop into deep folds. Instead, the smooth regular wrinkles transition into a temporary irregular shape (see the second empty circle), regular wrinkles again (see the filled circle), and finally a series of shapes involving non-smooth creases. In case 3 with $\mu_{ng}/\mu_{g}=1$ and $K_{ng}/K_{g}=1$ (right), no smooth wavy solutions are found when the circular configuration becomes locally unstable. When the circular configuration is still locally stable, the shape development that favors the lowest energy is constituted by a sequence of irregular non-smooth shapes before a regular non-smooth shape emerges.}\label{Fig_12}
\end{figure}
\begin{figure}[!b]
	\centering
	\includegraphics[width=1.0\textwidth]{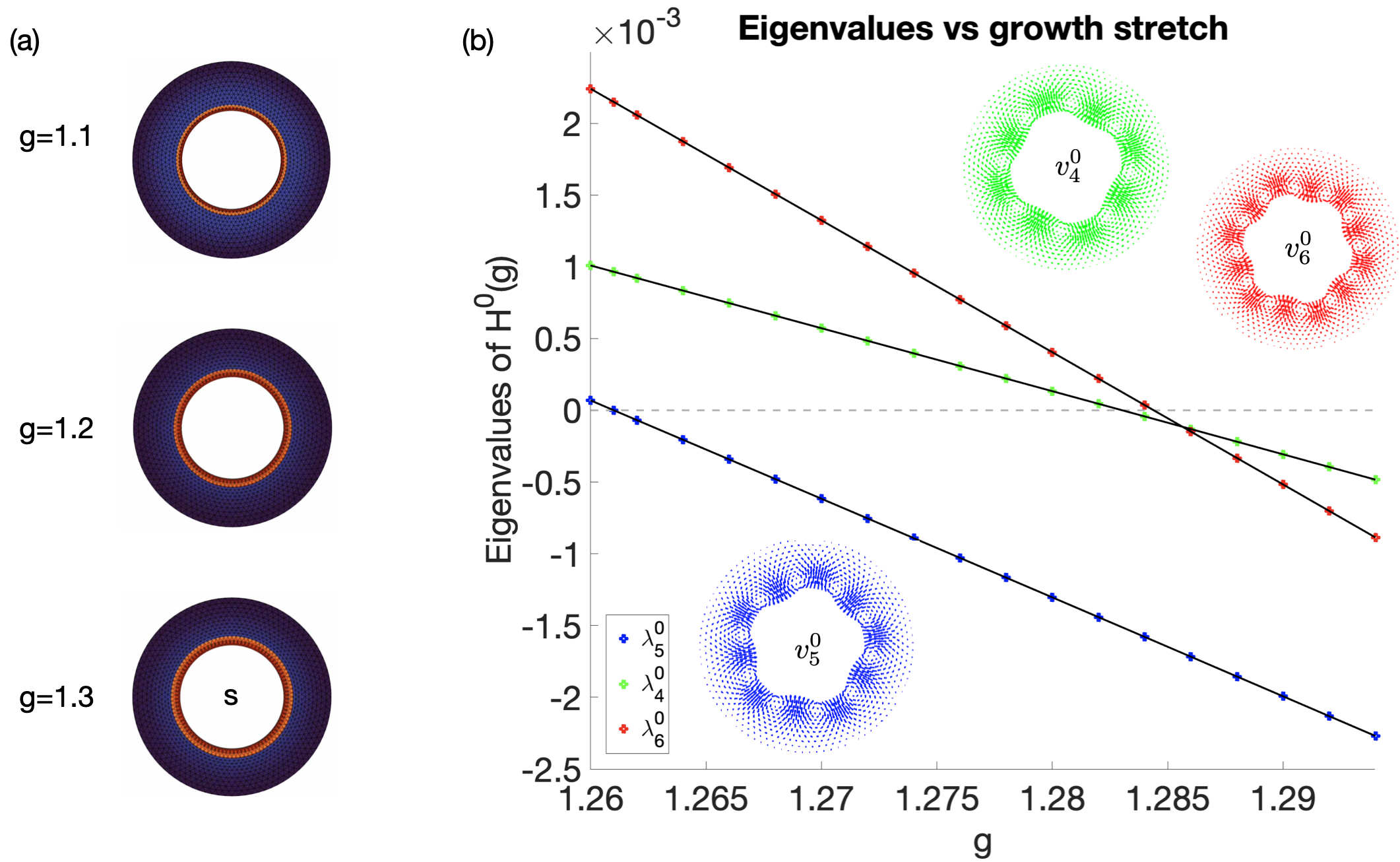}
	\caption{(a) Three snapshots of the annular configuration $s^0$ with different values of $g$ with $\mu_{ng}/\mu_{g}=0.1$ and $K_{ng}/K_{g}=0.1$. The placement of an ``s" at the center of a configuration indicates it is a saddle point at a specific value of $g$. (b) The three smallest eigenvalues of the Hessian $H(s^0)$  versus $g$.  The vector fields are generated by plotting the eigenvector components on the corresponding nodes. See text for details.}\label{Fig3}
\end{figure}

We apply Algorithm 1 to study shape development during growth in the problem described in Sec.\ref{sec:model}. In particular, we consider an annular region with its initial inner radius half of its outer radius, $R_{in} = 0.5$ and $R_{out}=1$, and triangulate it by $N_x=12$ layers of triangles along the radius, and $N_y = 92$ layers of triangles along the circumference. These triangles are similar, and increase their size from the interior boundary to the exterior boundary by a ratio of $(R_{out}/R_{in})^{1/N_x}\approx 1.06$ between successive layers. We set the growing domain $\Omega_0^{g}$ to be the first two interior layer of triangles (see Fig.\ref{Fig3}a) and the non-growing domain  $\Omega_0^{ng}$ in the rest exterior layers.  That is, we have $\mathbf{G} = g\mathbf{I}_2$ in the triangle elements of the first two layers, and $\mathbf{G} \equiv \mathbf{I}_2$ in all the other triangles exterior to the first two layers. We note that all the nodes satisfy the same local force Eq.(\ref{force}) and Hessian-component Eq.(\ref{Hessian}), and there is no special treatment on the nodes along the interface between the growing and non-growing domains.
 
Here, we list the notations that are frequently used in the upcoming results. The baseline annular configuration is denoted by $s^0$. The configurations that branches from $s^0$ are denoted by $s^m$, where $m$ is a wavenumber.  We use $s_m$ to denote configurations that have $m$ creasing indentations with regular spacings, and use $s_m^\Box$ to denote configurations that have $m$ creasing indentations with irregular  spacings.  The superscript $\Box$ is used to differentiate different irregular  solutions with the same number of creasing indentations.  For each configuration $s_m^\Box$, $\boldsymbol\varphi_{s_m^\Box}$ denotes its vectorized nodal coordinates, $H({s_m^\Box})$ denotes its Hessian, and $E(s_m^\Box)$ denotes its total elastic energy.

\subsection{Case 1: Formation of smooth wrinkles and foldings}
\label{case_1}

We first study the case where the elastic moduli in the growing region is much larger than the non-growing region. Following \cite{papastavrou2013mechanics}, we choose $\mu_{ng}/\mu_{g}=0.1$ and $K_{ng}/K_{g}=0.1$.  In this case, we find regular smooth wavy wrinkles branch from the circular configuration when it becomes locally unstable, which further develop into deep folds. See Fig.\ref{Fig_12} (left) for a schematic summary. Below we provide more details. 
\subsubsection{The first smooth-wrinkle configuration}
\begin{figure}[!b]
	\centering
	\includegraphics[width=0.7\textwidth]{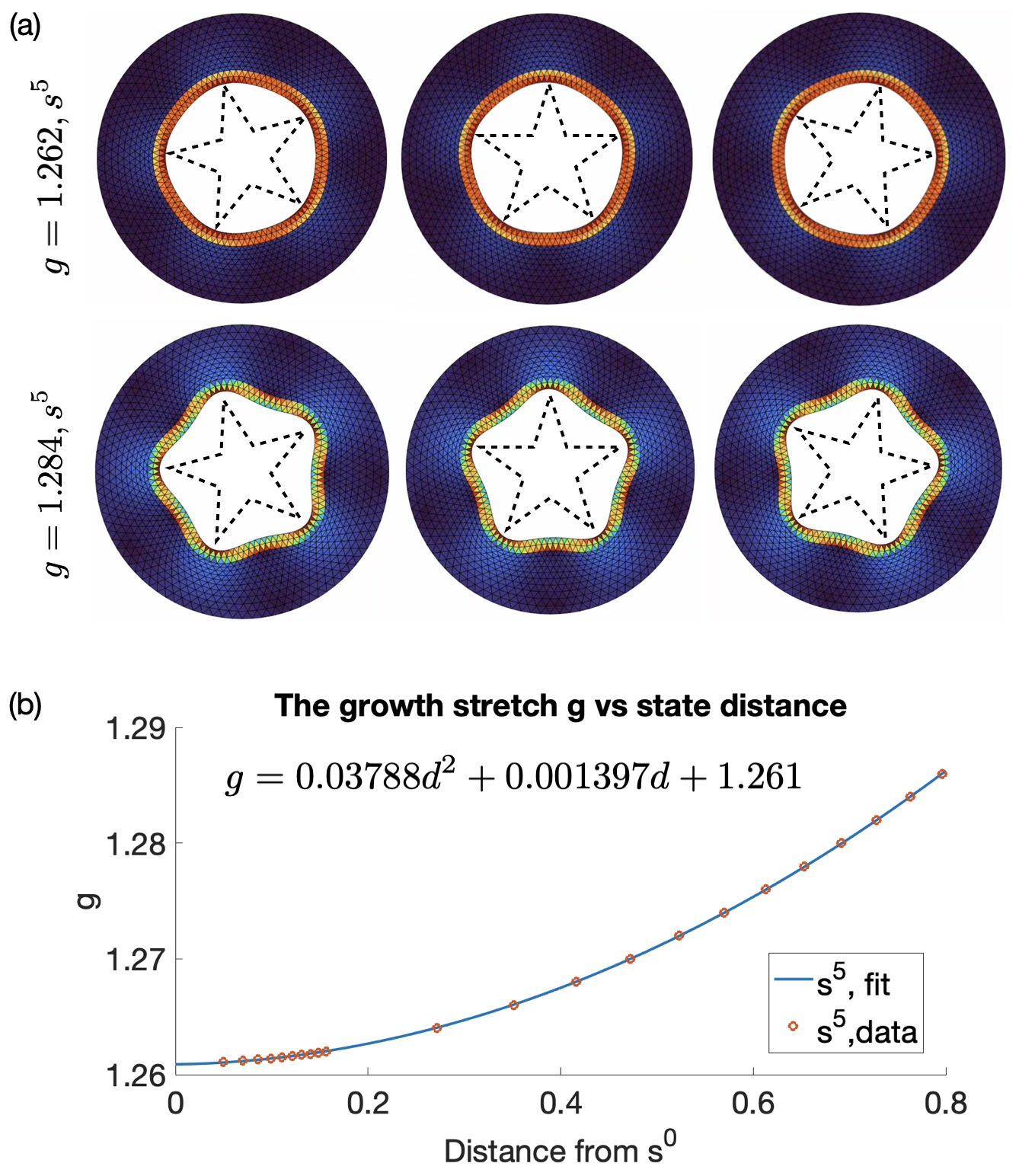}
	\caption{(a) Three isoenergetic configurations from $\{s^5\}$ at $g=1.262$ and $g=1.284$ with $\mu_{ng}/\mu_{g}=0.1$ and $K_{ng}/K_{g}=0.1$. The dashed stars with their vertices pointing to the indentations are added to help visualizing the difference in the phase of interior undulations. (b) The growth stretch parameter $g$ versus $d=||\boldsymbol{\varphi}_{s^5}-\boldsymbol{\varphi}_{s^0}||_2$, the distance between the configuration $\{s^5\}$ and $s^0$. \label{Fig4}}
\end{figure}

By increasing $g$ in the interior growing $\Omega_{g}$ from $g=1$, the annular configuration, denoted by $s^0$, continues to thicken in $\Omega_{g}$ (see Fig.\ref{Fig3}a). Its nodal-coordinate vector  $\boldsymbol\varphi_{s^0}$ has a positive-definite Hessian $H(s^0)$ until $g=1.261$ (Fig.\ref{Fig3}b). We find the eigenspace $V_5^0$ corresponding to  the smallest eigenvalue of $H(s^0)$ is a two-dimensional plane, and any vector $v_5^0\in V_5^0$ delineates a vector field of a wavenumber $m=5$ around the circumference of the annular domain (see example from the blue vector field in Fig.\ref{Fig3}b). Thus, $V_5^0$ corresponds to the first unstable wave mode from perturbing the annular configuration $s^0$. The fact that $V_5^0$ is two-dimensional is consistent with the expectation that there should be a family of shallow-wrinkled configurations with a wavenumber $m=5$ with only a shift in the phase. Indeed, by sampling three different eigenvectors $v_5^0\in V_5^0$ in the eigenplane, and perturbing $\boldsymbol\varphi_{s^0}\to\boldsymbol\varphi_{s^0}+\gamma v_5^0/|v_5^0|$ in the directions of these eigenvectors (with the magnitude $\gamma$), we obtained three new local energy minimizers (Fig.\ref{Fig4}a) with the same shape but differ in their phase of interior undulations. As the configurations correspond to undulations with $m=5$, we denote the set of configurations as $\{s^5\}$. 

By comparing the two rows in Fig.\ref{Fig4}a, it is obvious that the amplitude of undulations increases as the growth parameter $g$ increases. In fact, we can quantify how the set of configurations $\{s^5\}$ ``branch'' from $s^0$ by plotting the relation between $g$ and $d=||\boldsymbol{\varphi}_{s^5}-\boldsymbol{\varphi}_{s^0}||_2$, the distance between $\{s^5\}$ and $s^0$. See Fig.\ref{Fig4}b.  By fitting the $g$ versus $d$ by quadratic polynomials, we find the first-order coefficient is negligible compared to the second-order coefficient (see Fig.\ref{Fig4}b), suggesting a relation $(g-1.261)\propto d^2$. This indicates that the family $\{s^5\}$ is born from a pitchfork-like supercritical bifurcation at around $g_c = 1.261$, and this relation $d$ versus $g$ is analogous to the relation $\epsilon \propto \sqrt{g-g_c}$ in the weakly nonlinear analysis where $\epsilon$ is the buckling amplitude \cite{audoly2000elasticity,amar2023creases}. 


\subsubsection{More smooth-wrinkle solutions}
\begin{figure}[!b]
	\centering
	\includegraphics[width=1.0\textwidth]{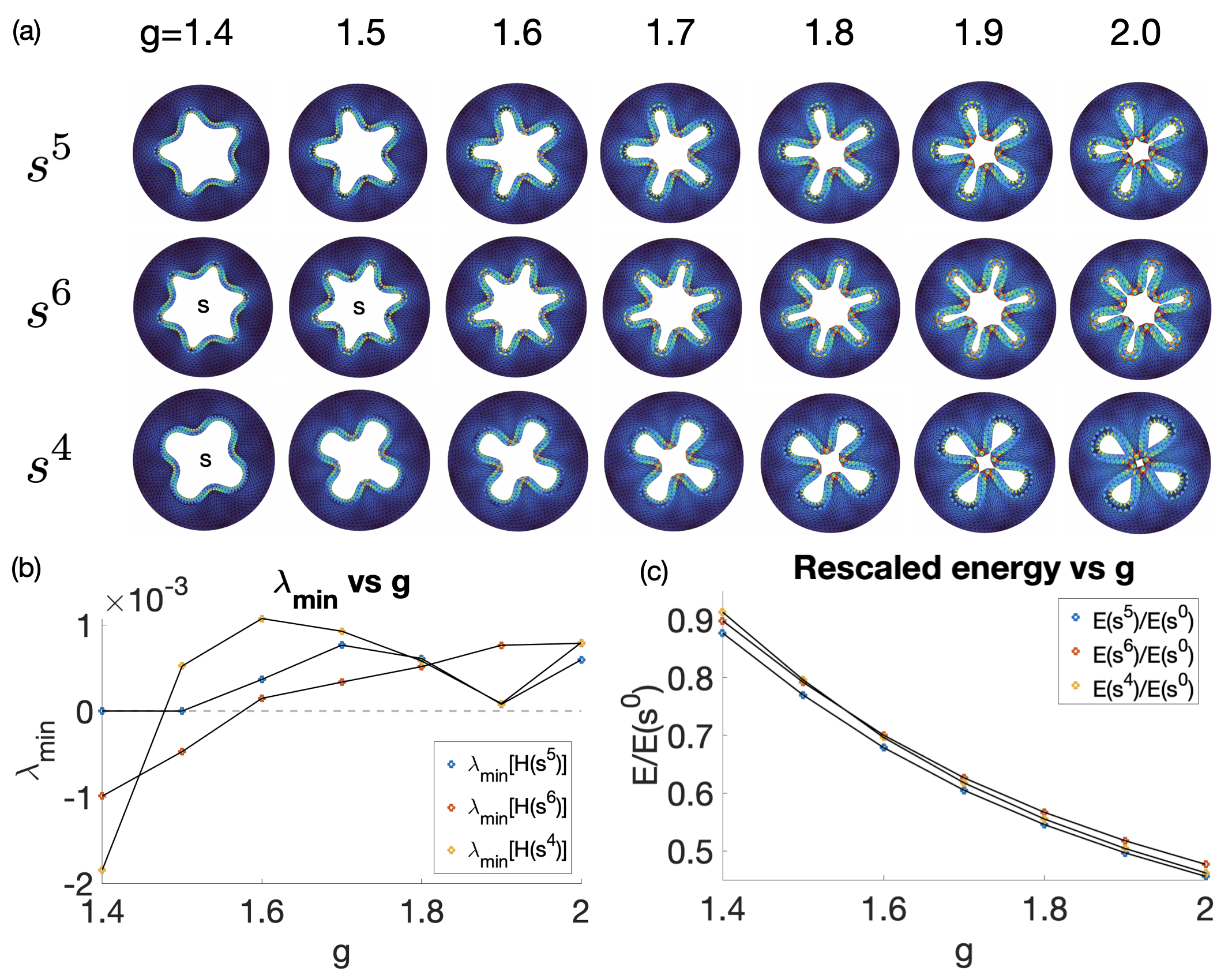}
	\caption{(a) The evolution of a $\{s^5\}$-,  $\{s^6\}$-, and $\{s^4\}$-configurations versus $g$ with $\mu_{ng}/\mu_{g}=0.1$ and $K_{ng}/K_{g}=0.1$. (b) The  smallest eigenvalues of $\{s^5\}$,  $\{s^6\}$, and $\{s^4\}$ versus $g$, respectively. (c) The energy of $\{s^5\}$,  $\{s^6\}$, and $\{s^4\}$ normalized by the energy $E(s^0)$ versus $g$. The placement of ``s" at the center of a configuration indicates it is a saddle point, not a local energy minimizer at the specific value of $g$. }\label{Fig6}
\end{figure}

Besides $\{s^5\}$, we find more configurations branching from $\{s^0\}$ at as $g$ further increases, namely $\{s^6\}$ and $\{s^4\}$. They are found by perturbing $\boldsymbol\varphi_{s^0}\to\boldsymbol\varphi_{s^0}+\gamma v_6^0/|v_6^0|$ and $\boldsymbol\varphi_{s^0}\to\boldsymbol\varphi_{s^0}+\gamma v_4^0/|v_4^0|$ when the second and third smallest eigenvalues of $H(s^0)$ become negative (see the green and red vector fields in Fig.\ref{Fig3}b). Similar wave patterns have been found previously in  \cite{papastavrou2013mechanics}. In addition to the previous results, we also find that the two new solutions $\{s^6\}$ and $\{s^4\}$ are both born as saddle points, which means they are not locally stable when they just branch from $s^0$  (Fig.\ref{Fig6}a and b). This is not surprising if the earlier born $\boldsymbol\varphi_{s^5}$ is close enough to the two new solutions  $\boldsymbol\varphi_{s^6}$ and $\boldsymbol\varphi_{s^4}$ with a lower energy. 

More interestingly, we find that $\{s^6\}$ and $\{s^4\}$ both become locally stable as $g$ further increase (Fig.\ref{Fig6}a and b). However, $\{s^5\}$ is still the absolute energy minimizer among the three family of solutions up to $g=2$ (Fig.\ref{Fig6}c). Strictly speaking, we cannot exclude the existence of other local minimizers that have a lower energy than $\{s^5\}$. We do not further increase $g$, as self-contact occur both in $\{s^5\}$ and $\{s^4\}$. We leave the consideration of treating self contact \cite{tallinen2016growth} in future work. 

\subsection{Case 2: Transition from regular wrinkles to creases, via temporary irregular  patterns}
\label{case_2}
Then, we adjust the ratio in elastic moduli between the growing and non-growing regions. In the direction of further decreasing  $\mu_{ng}/\mu_{g}$ and $K_{ng}/K_{g}$, we find similar patterns of smooth regular wrinkle patterns that later develop into deeper folds. This makes sense intuitively since the much-softer non-growing exterior layer is compliant enough to be distorted, without which the dramatic folds cannot be formed (Fig.\ref{Fig6}). However, in the direction of increasing $\mu_{ng}/\mu_{g}$ and $K_{ng}/K_{g}$, which means the exterior layer becomes less compliant, we find the smooth wrinkle patterns cannot develop into deeper folds as $g$ increases. In the case of  $\mu_{ng}/\mu_{g}=0.4$ and $K_{ng}/K_{g}=0.4$, we find successive transitory irregular patterns in the shape development of the growing annulus. See schematic summary in Fig.\ref{Fig_12}(middle) for the more complicated morphological evolution compared to case 1, and we show  details below.
\begin{figure}[!b]
	\centering
	\includegraphics[width=1.0\textwidth]{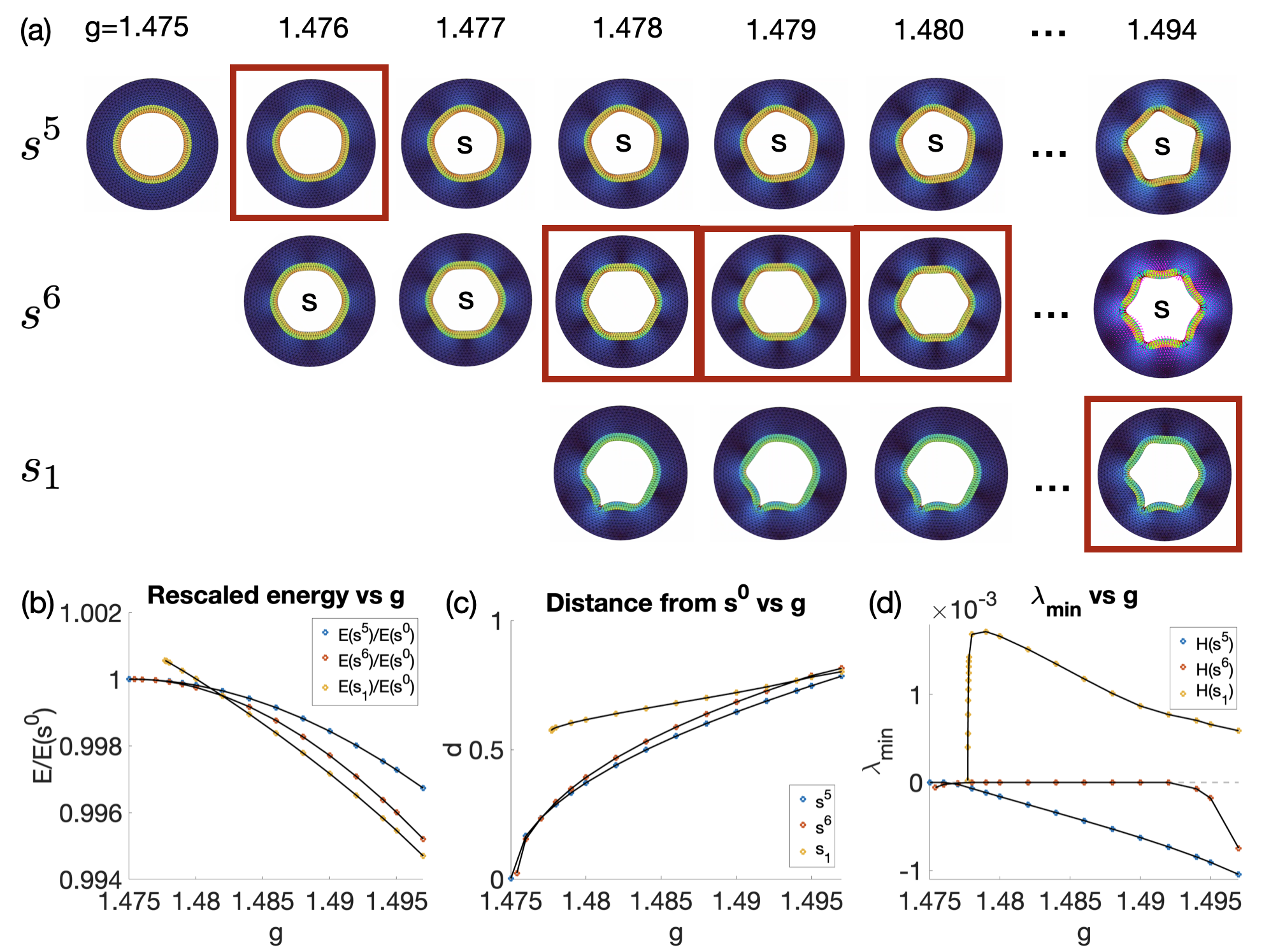}
	\caption{(a) The evolution of a $\{s^5\}$-,  $\{s^6\}$-, and $\{s_1\}$-configurations versus $g$ with $\mu_{ng}/\mu_{g}=0.4$ and $K_{ng}/K_{g}=0.4$. Red boxes highlight the solution with the smallest energy at each value of $g$. The placement of ``s" at the center of a configuration indicates it is a saddle point at a specific value of $g$. The purple vector field on $s^6(g=1.494)$ is from the eigenvector associated with the negative eigenvalue of $H(s^6)$.  (b) The energy of $\{s^5\}$,  $\{s^6\}$, and $\{s_1\}$ normalized by $E(s^0)$ versus $g$. (c) The distance of $\{s^5\}$,  $\{s^6\}$, and $\{s_1\}$ from the annular solution $s^0$: $||\boldsymbol\varphi_{s^5}-\boldsymbol\varphi_{s^0}||_2$, $||\boldsymbol\varphi_{s^6}-\boldsymbol\varphi_{s^0}||_2$, and $||\boldsymbol\varphi_{s_1}-\boldsymbol\varphi_{s^0}||_2$ versus $g$. (d) The smallest eigenvalues of $H(s^5)$, $H(s^6)$, and $H(s_1)$, respectively, versus $g$.}\label{Fig7}
\end{figure}
\subsubsection{A smooth transitory irregular  solution}
First, we find a transitory irregular  smooth solution that has not been reported before. Similar to case 1 in Sec.\ref{case_1}, we find the first new smooth solution $\{s^5\}$ branching from the annular solution $s^0$ (see first row in Fig.\ref{Fig7}a). Following that, we find the second new solution $\{s^6\}$ branching from $s^0$ (see second row in Fig.\ref{Fig7}a). See also Fig.\ref{Fig7}b and c as confirmation that the two solution states $\{s^5\}$ and $\{s^6\}$ indeed branch out from $s^0$. 
In contrast to Sec.\ref{case_1}, $\{s^5\}$ becomes locally unstable while $\{s^6\}$ becomes locally stable in the interval of $g\in[1.476 1.478]$. Within that interval, we find a transitory irregular  configuration $\{s_*\}$ that branches from $\{s^5\}$ exactly when $\{s^5\}$ becomes a saddle point while approaching the saddle point $\{s^6\}$. It merges with $\{s^6\}$ exactly when $\{s^6\}$ becomes a local minimizer. See Fig.\ref{Fig7}a for the shape dynamics of $\{s_*\}$, Fig.\ref{Fig7}b and d for the confirmation of its local stability when $\{s^5\}$ and $\{s^6\}$ are both saddle points, and Fig. \ref{Fig7}c for its emergence from $\{s^5\}$ and the later merging with $\{s^6\}$.

To the best of our knowledge, this is the first irregular  smooth pattern found from growth-elasticity modeling. We argue the failure to find such a solution previously is due to the limitation of Newton-iteration-based numerical methods \cite{papastavrou2013mechanics}.  Since $\{s^5\}$ or $\{s^6\}$ are very close already, using Newton's method alone fails to identify $\{s_*\}$--when we started with one solution, say $s^5$, we end up finding the other, $s^6$.  The energy-stable gradient-descent steps are necessary to localize the search to the energy valley between $\{s^5\}$ and $\{s^6\}$. 
\begin{figure}[!b]
	\centering
	\includegraphics[width=1.0\textwidth]{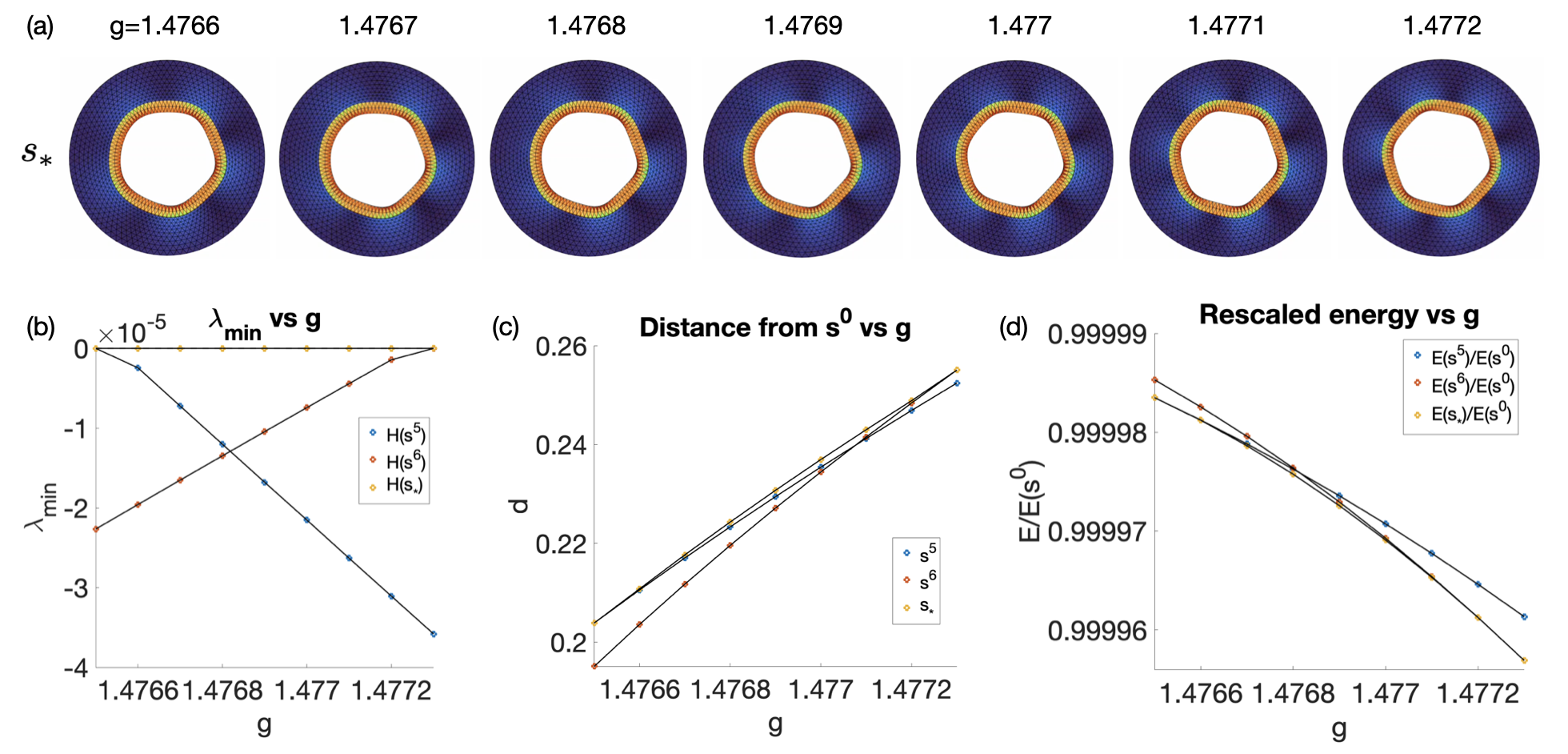}
	\caption{(a) The evolution of a $\{s_*\}$-configuration versus $g$, incrementally transitioning from a  $\{s^5\}$-configuration to a $\{s^6\}$-configuration with $\mu_{ng}/\mu_{g}=0.4$ and $K_{ng}/K_{g}=0.4$. (b) The smallest eigenvalue of $H(s_*)$ versus $g$, juxtaposed with those of $H(s^5)$ and $H(s^6)$. This shows $\{s_*\}$ is the only locally stable solution while the other two are saddle points. (c) The distance of $\{s_*\}$ from the annular solution $s^0$: $||\boldsymbol\varphi_{s_*}-\boldsymbol\varphi_{s^0}||_2$, juxtaposed with those of $\{s^5\}$ to $\{s^6\}$. Notice $||\boldsymbol\varphi_{s_*}-\boldsymbol\varphi_{s^0}||_2$ converges to $||\boldsymbol\varphi_{s^5}-\boldsymbol\varphi_{s^0}||_2$ at the lower end of $g$, and converges to $||\boldsymbol\varphi_{s^6}-\boldsymbol\varphi_{s^0}||_2$ at the higher end. (d) The normalized energy $E(s_*)/E(s^0)$ versus $g$, juxtaposed with those of $E(s^5)$ and $E(s^6)$. Notice $\{s_*\}$ has the lowest energy throughout its existence range.}\label{Fig8}
\end{figure}

\subsubsection{Transient irregular solutions with a mixture of smooth and non-smooth indentations}
Both $\{s^5\}$ and $\{s^6\}$ become unstable as $g$ increases (see the first two rows when $g=1.494$). This makes intuitive sense as the external layer becomes less compliant, forbidding deeper folds.  Instead, we find an irregular  pattern as the new local energy minimizer $\{s_1\}$ (see the third row from Fig.\ref{Fig7}a). The solution $\{s_1\}$ is found by perturbing the saddle point $\{s^6\}$ at $g=1.494$ by the vector field from the unstable eigenspace (see the purple vector field in Fig.\ref{Fig7}a) followed by Algorithm 1. It is essentially different from the previous solutions $\{s^5\}$ and $\{s^6\}$ such that it has a non-smooth indentation (crease) and five smooth indentation along the inner boundary. 

Interestingly, $\{s_1\}$ is not branched from $\{s^6\}$, or any other critical points. We can see this by decreasing the value of $g$ on $\{s_1\}$ continuously, and find it is born earlier around $g=1.478$. See the third row from Fig.\ref{Fig7}a. We can further confirm that it is born from a distortion of the energy landscape at somewhere away from any of the pre-existing critical points $\{s^0\}$, $\{s^5\}$, and $\{s^6\}$, by checking its energy dynamics (Fig.\ref{Fig7}b), its distance from $s^0$ in comparison with other solutions (Fig.\ref{Fig7}c), and the smallest eigenvalue of its Hessian $H(s_1)$) (Fig.\ref{Fig7}d)-- in fact it is born with a higher energy value than any of other configurations, and gradually becomes more energetically favorable as $g$ increases; at birth, its distance from $\{s^0\}$ is finite and there is a sharp increase of $H(s_1)$'s smallest eigenvalue from zero, which is again qualitatively different from those in $\{s^5\}$ and $\{s^6\}$. 

 \begin{figure}[!b]
	\centering
	\includegraphics[width=1.0\textwidth]{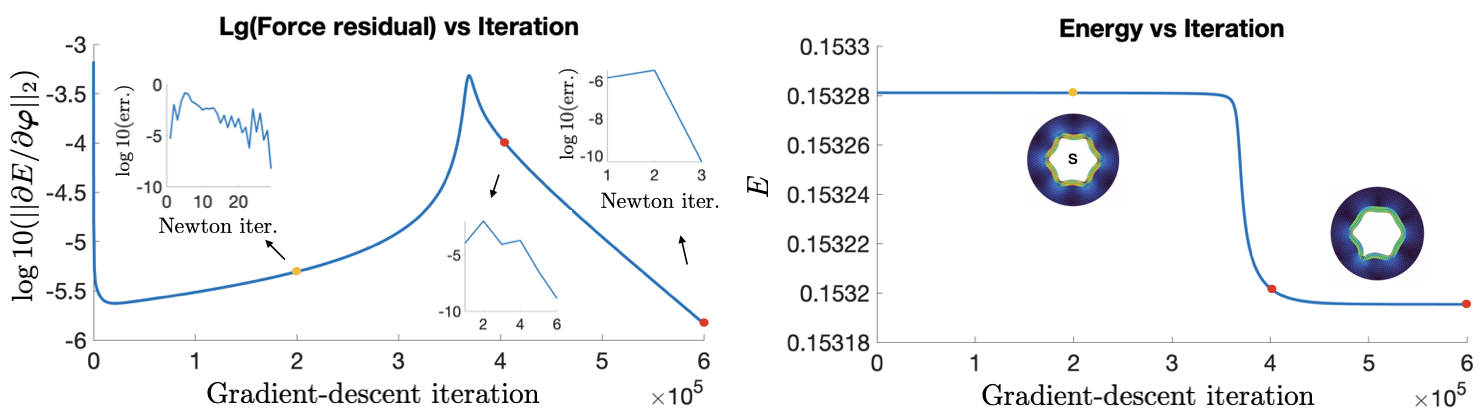}
	\caption{The decimal logarithm of the total force residual (left) and energy (right) versus the gradient-descent iteration. As described in Algorithm 1, we perform the gradient-descent iteration and then perform Newton's iteration upon a chosen $MaxIter_1$. Yellow and red dots tag when trials of Newton's iteration are performed. When Newton's iteration is introduced too early ($MaxIter_1=2\times 10^5$), the inset shows the iterates approach the saddle solution $\{s^6\}$ instead of the new local minimizer $\{s_1\}$. When $MaxIter_1=4\times 10^5$ and $6\times 10^5$, the Newton's iterates approach $\{s_1\}$ with force residual $||\partial E/\partial \boldsymbol{\varphi} ||_2<10^{-7}$ within $6$ and $3$ Newton's iteration, respectively. }\label{Fig_revision}
\end{figure}
We emphasize that the gradient-descent iteration from Algorithm 1 is essential to find $\{s_1\}$. In Fig. \ref{Fig_revision}, we show the global force-residual $||\partial E/\partial \boldsymbol{\varphi} ||_2$ and total energy $E$ versus the iteration. The gradient-descent iteration is critical to guide the numerical iterates to deviate from $\{s^6\}$ to approach $\{s_1\}$. If we perform the Newton's iteration too early (e.g., $Iter.=2\times 10^5$), the Newton's iterates will go back to the saddle point $\{s^6\}$. Only when we perform Newton's iteration after the gradient-descent iterates passes a steep region of the energy landscape (e.g., $Iter.=4\times 10^5 \text{ or }6\times 10^5$), the Newton's iterates can approach $\{s_1\}$. Thus, previous Newton-iteration-based numerical methods \cite{papastavrou2013mechanics} are not able to find such a solution, which might be rescued by pseudo-arclength approaches \cite{groh2022morphoelastic}. However, the pseudo-arclength method is much more expensive than our method, and may fail to resolve the situations when multiple solution paths are too close (e.g., $\{s_*\}$ versus $\{s^6\}$ and $\{s^5\}$). The usage of energy landscape is essential to identify different configurations robustly. 


 Mathematically, the crease is very different from the smooth wrinkles. Focusing on the material points along the interior boundary, the crease formation involves a transformation from a smooth curve (starts with a circle) to an non-smooth curve with the local tangent $\hat{t}$ not defined at the crease. This implies that $\nabla_0\varphi$ is singular at the crease since $\hat{t} = \hat{t}_0\cdot \nabla_0\varphi/|\hat{t}_0\cdot \nabla_0\varphi|$, where $\hat{t}_0$ is the local tangent along the initial circle. This is not the case for smooth wrinkles. We note that having singularity in $\nabla_0 \varphi$ is not an issue in an energetic formulation-- both for Galerkin and Rayleigh-Ritz formulations as long as the singularity is taking a zero measure. 

\begin{figure}[!b]
	\centering
	\includegraphics[width=1.0\textwidth]{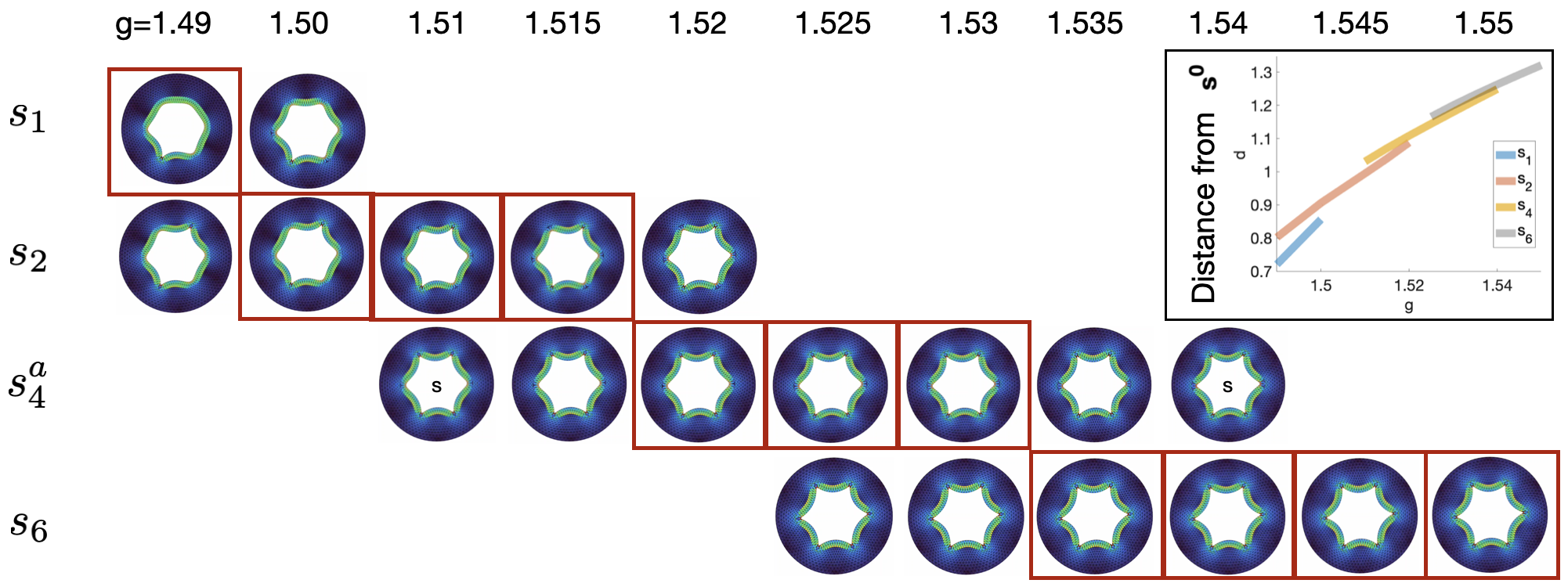}
	\caption{The evolution of configurations with more creases from $\{s_1\}\to\{s_2\}\to\{s_4^a\}\to\{s_6\}$ with $\mu_{ng}/\mu_{g}=0.4$ and $K_{ng}/K_{g}=0.4$. Red boxes highlight the configuration with the smallest energy at each value of $g$. The placement of ``s" at the center of a configuration indicates it is a saddle point at a specific value of $g$. The inset (black box) shows the distance of these configuration from the annular solution versus $g$.}\label{Fig9}
\end{figure}

 As we continue to increase $g$, we find more irregular  solutions with a mixture of sharp creases and smooth indentations. Here we only show $\{s_2\}$ with two opposing creases and four smooth indentations and $\{s_4^a\}$ with two opposing smooth indentation and four creases. See the second and third rows of Fig.\ref{Fig9}. They both only live in short interval of $g$. As $g$ further increase, a regular solution $\{s_6\}$ with six creases are found. 
 
 One can see from Fig.\ref{Fig9} that there exist at least two solutions for the same value of $g$,  among which we generate a minimal-energy path (solutions in red boxes). In fact, we find more local minimizers such as $s^5(g=1.545)$ and $s^7(g=1.55)$, and a saddle point $s^4(g=1.55)$ as well as other irregular  patterns with creases, but none of them has a lower energy than the minimal energy among $\{s_1\}$, $\{s_2\}$, $\{s_4^a\}$, and $\{s_6\}$ at a specific value of $g$ when they are available. 
 
From the configurations shown in Figs.\ref{Fig6} and \ref{Fig9}, our result suggests a minimal-energy path of shape development:  $\{s^5\}\to\{s_*\}\to\{s^6\}\to\{s_1\}\to\{s_2\}\to\{s^a_4\}\to\{s_6\}$. To the best of our knowledge, previous simulations have not captured transitory patterns with a mixture of smooth and non-smooth indentations, while it is observed in reality. See Fig. \ref{Fig-1} for example. As shown in Fig.\ref{Fig8}, the usage of the energy formulation is essential to the finding of these irregular  solutions.

\subsection{Case 3: Transitions from irregular  creasing patterns to regular  ones}
\label{case_3}
In the previous two cases (see summary in Fig.\ref{Fig_12} left and middle), the initial shape transition from the circular one involves the regular shallow wrinkles. In this subsection, we show that regular shallow wrinkles do not have to be the initial symmetry-breaking configurations in the case of equal stiffness among the two layers  ($\mu_{ng}/\mu_{g}=1$ and $K_{ng}/K_{g}=1$). In particular, we find a large number of  creasing patterns which altogether provide a minimal-energy shape evolution by increasing the number of creases. See summary in Fig.\ref{Fig_12} (right) and we provide more details below.
\subsubsection{Multiple regular creasing patterns}

\begin{figure}[!b]
	\centering
	\includegraphics[width=1.0\textwidth]{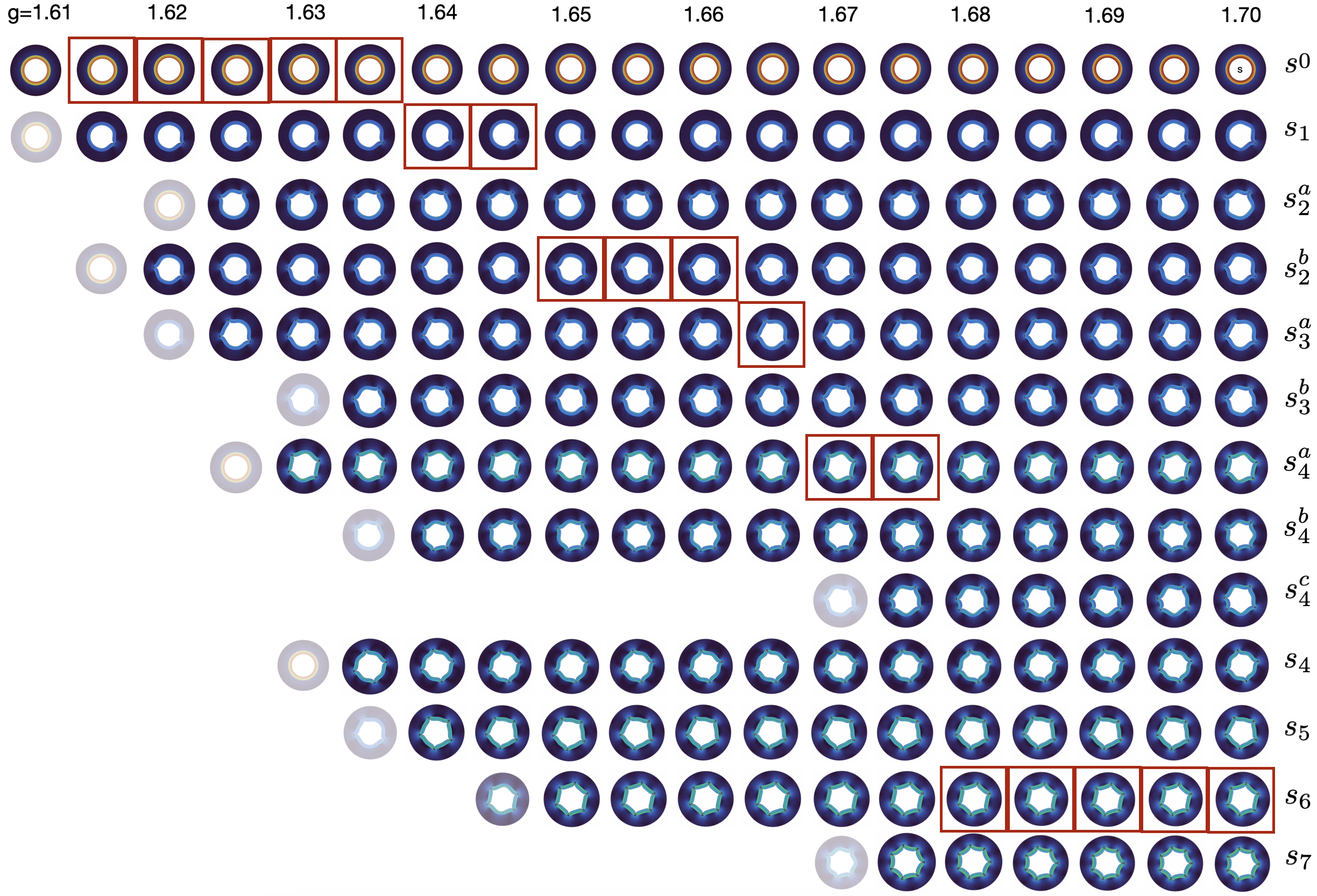}
	\caption{The evolution of different crease formations with $\mu_{ng}/\mu_{g}=1$ and $K_{ng}/K_{g}=1$. Each row corresponds to one local-minimizer configuration and its coverage along the change of $g$, and each row is preceded by a more translucent configuration indicating another branch of state that is found when the current solution disappears as $g$ further decreases. Red boxes highlight the shape with the smallest energy at each value of $g$. All solutions are local energy minimizers except the annular solution $s^0$ at $g=1.7$ with a ``s" at the center of the configuration. 
}\label{Fig10}
\end{figure}

First of all, we find the annular solution $s^0$ is locally stable up to when $g=1.7$ (see Fig.\ref{Fig10}, top row $s^0$). The 2D eigenspace of the Hessian $H(s^0)$ corresponding to the first negative eigenvalue presents us a family of vector fields with a wavenumber $m=6$, similar to the vector field $v_6^0$ shown in Fig.\ref{Fig3}b. This again suggests a family of new solutions $\{s^6\}$ that provide the same shape pattern but differ by a rotation. However, when we perturb $\boldsymbol\varphi_{s^0}\to\boldsymbol\varphi_{s^0}+\gamma v_6^0/|v_6^0|$ with a range of small magnitudes $\gamma=10^{-10}-10^{-3}$, we have not found a smooth wrinkled solution $\{s^6\}$ close to $s^0$. Instead, we always find the creasing configuration $\{s_6\}$ with a finite distance from $s^0$. See Fig.\ref{Fig10}, the last second row $s_6$ at $g=1.7$, for a configuration $s_6$. 

\begin{figure}[!b]
	\centering
	\includegraphics[width=1.0\textwidth]{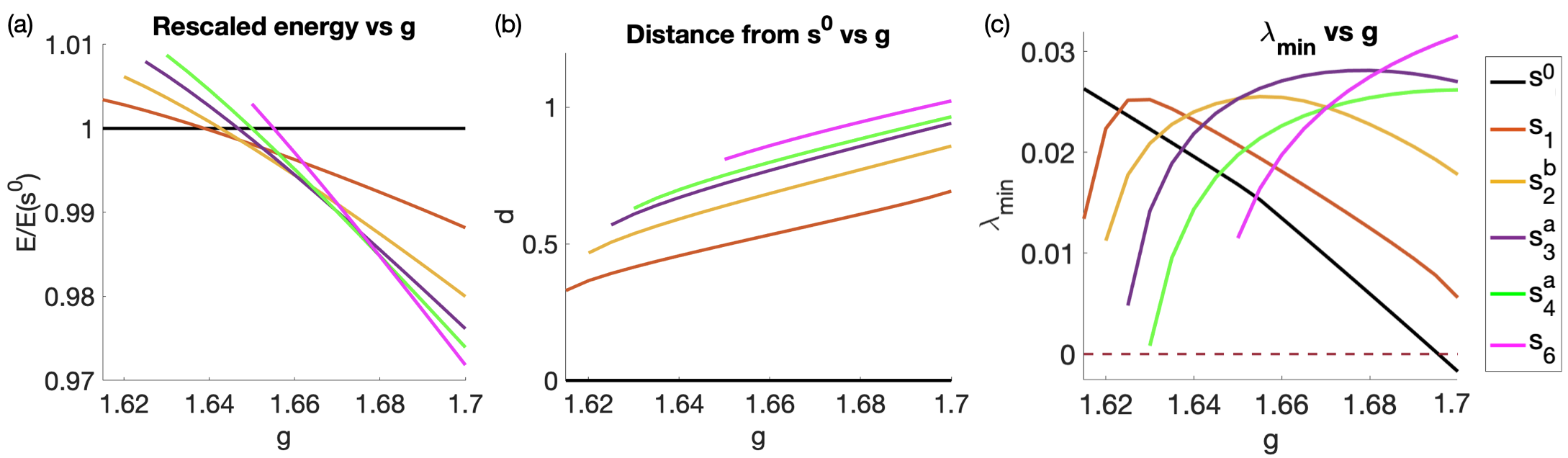}
	\caption{(a) The normalized energy $E/E(s^0)$ versus $g$ among $s^0$, $\{s_1\}$, $\{s_2^b\}$, $\{s_3^a\}$, $\{s_4^a\}$, and $\{s_6\}$. (b) The distance of $s^0$, $\{s_1\}$, $\{s_2^b\}$, $\{s_3^a\}$, $\{s_4^a\}$, and $\{s_6\}$ from $s^0$.  (c) The smallest eigenvalue of $H(s^0)$, $H(s_1)$, $H(s_2^b)$, $H(s_3^a)$, $H(s_4^a)$, and $H(s_6)$ versus $g$.}\label{Fig11}
\end{figure}

As we have seen from Sec.\ref{case_2}, this creasing configuration may be born earlier. Thus we decrease $g$ continuously on $\{s_6\}$, and find $\{s_6\}$ occurs as early as $g=1.65$ (see Fig.\ref{Fig10}, row $s_6$), with a higher energy than $s^0$ and a finite distance from $s^0$ (Fig.\ref{Fig11}). Thus we conclude that when $s^0$ becomes unstable at $g=1.7$, there is no new solution branching from it. Rather, the other local energy minimizer $\{s_6\}$ already exists and the energy barrier between $s^0$ and $\{s_6\}$ is removed when $s^0$ becomes unstable. As the route from $s^0$ to $\{s_6\}$  is aligned with the unstable eigenvector $v_6^0$, we find $\{s^6\}$ by perturbing $s^0$ along $v_6^0$.

As there is no new solutions branching from $s^0$ at $g=1.7$ when it becomes unstable, we suspect there are other local energy minimizers existing at $g=1.7$ (and even earlier), but are divided from $s^0$ by energy barriers. Thus, we assume the directions along eigenvectors with smaller eigenvalues (though positive) corresponds to lower energy barriers between $s^0$ and another local minimizer. Thus, we perturb $\boldsymbol\varphi_{s^0}\to\boldsymbol\varphi_{s^0}+\gamma v_x^0/|v_x^0|$ along eigenvectors $v_x^0$'s with the ascending order of the eigenvalues with a finite magnitude $\gamma$. We choose $\gamma =1$ such that the perturbation is large enough to pass some energy barriers but will not generate unphysical configurations. Along the direction corresponding to the 2nd-smallest (yet positive) eigenvalue, we find another regular configuration $\{s_5\}$ at $g=1.7$. Similarly, we find regular configurations $\{s_7\}$ and $\{s_4\}$ by perturbing $s^0$ along directions with the 3rd- and 4th-smallest eigenvalue. See row $s_4$, $s_5$, and $s_7$ in Fig.\ref{Fig10} for their configurations. 

We have not found more regular configurations by keep using the eigenvectors of larger eigenvalues. For example, by perturbing $s^0$ along directions with the 5th-smallest eigenvalue, we still find $\{s_4\}$, and along the direction of the 6th-smallest eigenvalue, we only find the saddle point $\{s^0\}$. Interestingly, along the direction of the 7th-smallest eigenvalue, we find an irregular  configuration $\{s_4^c\}$, which still has four indentations, but with irregular spacings in between. See row $s_4^c$ in Fig.\ref{Fig10}.

\subsubsection{More irregular  creasing patterns}
More irregular  configurations are found by decreasing the value of $g$ with a decrement $\Delta g=-0.005$ on the configurations $\{s_4^c\}$, $\{s_5\}$, $\{s_6\}$, and $\{s_7\}$ until they disappear. Respectively, they lead to the irregular  solutions $\{s_3^b\}$, $\{s_3^a\}$, $\{s_4^a\}$, and $\{s_4^b\}$. We note that the disappearance of $\{s_4\}$ at $g=1.63$ only leads to the annular solution $s^0$. We find more irregular  solutions $\{s_2^b\}$, $\{s_1\}$, and $\{s_2^a\}$ by decreasing the value of $g$ on $\{s_3^b\}$, $\{s_3^a\}$, and $\{s_4^b\}$, respectively. See Fig.\ref{Fig10} and its caption for details.

All the irregular  configurations are born initially with a higher energy than $s^0$ and with a finite distance from $s^0$. See Fig.\ref{Fig11} for the details of $\{s_1\}$, $\{s_2^b\}$, $\{s_3^a\}$, $\{s_4^a\}$ for example. Among all these solution branches, $\{s_1\}$ with a singular crease indentation appears the earliest at $g=1.615$, and its energy $E(s_1)$ becomes more favorable (lower) than $E(s^0)$ the earliest at $g=1.64$. This means that the very first symmetry breaking we have found is a localized indentation $\{s_1\}$, instead of regular wrinkles or creasing patterns.

As $g$ increases, other solution branches become more energetically-favorable, leading to a configuration path  $\{s_1\}\to\{s_2^b\}\to\{s_3^a\}\to\{s_4^a\}\to\{s_6\}$. Along this path, the first four, as $g$ increases, are all irregular. The regular  $\{s_6\}$ is only selected later at $g=1.68$. This situation is different from the extreme setting where the external layer is rigid, and only regular crease patterns are found theoretically \citep{jin2011creases}. This situation is also very different from the cases in Sec.\ref{case_1} and Sec.\ref{case_2}. Compared to Sec.\ref{case_1}, the essential difference is that there is no smooth undulating solution found. Compared to Sec.\ref{case_2}, there are two important differences. First, the initial symmetry breaking leads to an irregular  configuration that resembles a localized mechanical effect. Second, the transitions among solution branches $\{s_2^b\}\to\{s_3^a\}\to\{s_4^a\}\to\{s_6\}$ seem random, as they involve not only the change of crease numbers but also the redistribution of spacing among the indentations. This is in contrast to Fig.\ref{Fig9} where the shape transitions only involve smooth indentation turning into creases while the spacing between the six indentations is well maintained. 

 Lastly, we emphasize that the annular solution $s^0$ remains stable until $g=1.7$ after a bountiful of non-smooth solutions are already born. See smallest-eigenvalue of $H(s^0)$ versus $g$ in Fig.\ref{Fig11}c.  This is different from the two cases in Sec.\ref{case_1} and Sec.\ref{case_2} where the first new solution is a wrinkled configuration. The wrinkled solution branches from the annular solution when it becomes locally unstable with negative Hessian eigenvalues as shown in Fig.\ref{Fig3}. Since non-smooth solutions are not branched from the annular solution, the positivity of the Hessian $H(s^0)$ does not dictate the emergence of new non-smooth solutions.

\section{Conclusion}
\label{conclusion}
In this work, we introduced an energy-optimization FEM method to study growth-elasticity models. We applied this method to simulate the growth-driven shape dynamics of a bilayer annulus, which resembles the 2D cross-section of gut tubes. Our results demonstrate that the new method is capable of capturing regular wrinkles and their nonlinear development into deep folds, similar to the findings of \cite{papastavrou2013mechanics}, as well as regular creasing patterns, akin to those reported by \cite{jin2011creases}. Beyond these previously reported configurations, we have identified new irregular  patterns featuring smooth wrinkles, a mixture of smooth wrinkles and non-smooth creases, and creases with non-uniform spacings.

These irregular  patterns are observable in the cross sections of mammalian gut tubes (see Fig.\ref{Fig-1} and more data in \citep{li2011surface}), but they have not been theoretically recovered in previous works. One might argue that the non-uniform spacing of the undulations and the mixture of smooth and non-smooth indentations are due to initial geometric or material defects, small variations in growth along the circumference, or other random forces. However, based on our new numerical results, we have shown that these irregular  configurations can be explained within the framework of growth-elasticity without complicating the circumferential growth patterning or adding other mechanical factors.

We argue that the primary reason these solutions were not recovered by previous methods lies in the Galerkin-type formulation followed by Newton's iteration. This type of method disregards the energy information that could be used to design more robust solution search strategies. Starting from the energetic formulation of growth-elasticity, our method is able to search for solutions more robustly. It can identify a new energy-minimizer that is very close to two saddle points (Fig.\ref{Fig8}) and also discover new energy-minimizers that are far from the initial guess (e.g., Fig.\ref{Fig_revision}). These scenarios are challenging for Galerkin-type formulations and Newton's iteration, which might be improved by implementing pseudo-arclength numerical continuation \citep{groh2022morphoelastic}.

To conclude, our energy-optimization method significantly expands the explanatory power of growth-elasticity, even with simple growth patterning. Theoretically, how do we interpret the simulated shape evolutions among multiple solutions, given that live-imaging crypt formation data are rare? One natural answer is to select the lowest-energy configuration at each time and connect them sequentially. Fig.\ref{Fig6} provides the simplest continuous scenario: a five-fold smooth configuration is selected as the energy-favorable shape, and the five folds deepen into smooth crypts as growth increases. The scenarios with irregular  snapshots are more complex. For example, the shape evolutions depicted in Figs.\ref{Fig9} and \ref{Fig10}. Fig.\ref{Fig9} illustrates a shape evolution where local smooth indentations transform into sharp ones, offering a more or less continuous dynamic process without significant disruption of successive configurations. This is not the case in Fig.\ref{Fig10}, where the number of creases increases while some of the pre-existing creases relocate. This results in a highly variable developmental process, which may be more applicable to abnormal intestinal growth rather than normal crypt formation.

\appendix
\section{The derivation of Eq.(\ref{force})}
We derive Eq.(\ref{force}) by the following steps:
\begin{eqnarray}
\nonumber
\frac{\partial E^{\alpha_k}}{\partial \varphi^k}&=& J_g^{\alpha_k}\big(  \frac{\partial W^{\alpha_k}}{\partial \mathbf{F}_e^{\alpha_k}}: \frac{\partial \mathbf{F}_e^{\alpha_{k}}}{\partial \varphi^k}\big)A_0^{\alpha_k}\\\nonumber
&=&J_g^{\alpha_k}\big(  \frac{\partial W^{\alpha_k}}{\partial \mathbf{F}_e^{\alpha_k}}: \frac{\partial \mathbf{m}^{\alpha_{k}}}{\partial \varphi^k}[\mathbf{M}^{\alpha_k}]^{-1}[\mathbf{G}^{\alpha_k}]^{-1} \big)A_0^{\alpha_k}\\\nonumber
&=&J_g^{\alpha_k}\big(  \underbrace{\frac{\partial W^{\alpha_k}}{\partial \mathbf{F}_e^{\alpha_k}} [\mathbf{F}_e^{\alpha_{k}}]^{T}}_{J_e^{\alpha_k}\boldsymbol{\sigma}^{\alpha_k} }:\frac{\partial \mathbf{m}^{\alpha_k}}{\partial \varphi^k} [\mathbf{m}^{\alpha_k}]^{-1}\big)A_0^{\alpha_k}\\\nonumber
&=&J_g^{\alpha_k}{J_e^{\alpha_k}A_0^{\alpha_k}\boldsymbol{\sigma}^{\alpha_k}} :\frac{\partial \mathbf{m}^{\alpha_k}}{\partial \varphi^k} [\mathbf{m}^{\alpha_k}]^{-1}\\\nonumber
&=&A^{\alpha_k}\boldsymbol{\sigma}^{\alpha_k} :\frac{\partial \mathbf{m}^{\alpha_k}}{\partial \varphi^k} [\mathbf{m}^{\alpha_k}]^{-1}\\\nonumber
&=&\frac{1}{2}\boldsymbol{\sigma}^{\alpha_k} :\big(\det{\mathbf{m}^{\alpha_k}}\frac{\partial \mathbf{m}^{\alpha_k}}{\partial \varphi^k} [\mathbf{m}^{\alpha_k}]^{-1}\big)\\\nonumber
  &=&\frac12\boldsymbol{\sigma}^{\alpha_k}\varphi^{k\perp,\alpha_k}.
\end{eqnarray}
In the third line, $J_e^{\alpha_k} =\det (\mathbf{F}_e^{\alpha_k})$, and we have used the Cauchy stress $\boldsymbol{\sigma}^{\alpha_k} =(J_e^{\alpha_k})^{-1}\frac{\partial W^{\alpha_k}}{\partial \mathbf{F}_e^{\alpha_k}} [\mathbf{F}_e^{\alpha_{k}}]^{T}$. In the fifth equality, we have used the fact that the area variation $A^{\alpha_k}/A_0^{\alpha_k} = J_g^{\alpha_k}J_e^{\alpha_k}$. 

\section{The derivation of Eq.(\ref{stable}) and Eq.(\ref{optimal})}
Here, we write $ \nabla E(\boldsymbol{\varphi}_{i}):=\partial E/\partial \boldsymbol{\varphi}_i$ and $ \nabla^2 E(\boldsymbol{\varphi}_{i}):=\partial^2 E/(\partial \boldsymbol{\varphi}_i\partial \boldsymbol{\varphi}_i)$. First we derive  $\Delta s_{stable}$ from Eq.(\ref{stable}) for the iteration
\begin{eqnarray}
\boldsymbol{\varphi}_{i+1}-\boldsymbol{\varphi}_{i} = -{\Delta s} \nabla E(\boldsymbol{\varphi}_{i}).
\label{iter_rule}
\end{eqnarray}
Assuming $||\boldsymbol{\varphi}_{i+1}-\boldsymbol{\varphi}_{i}||_2\ll 1$, we see from the Taylor expansion 
\begin{equation}
\nonumber
E(\boldsymbol{\varphi}_{i+1})-E(\boldsymbol{\varphi}_{i}) = \nabla E(\boldsymbol{\varphi}_{i})(\boldsymbol{\varphi}_{i+1}-\boldsymbol{\varphi}_{i})+\frac{1}{2}(\boldsymbol{\varphi}_{i+1}-\boldsymbol{\varphi}_{i})^T\nabla^2 E(\boldsymbol{\varphi}_{i})(\boldsymbol{\varphi}_{i+1}-\boldsymbol{\varphi}_{i})+o(||\boldsymbol{\varphi}_{i+1}-\boldsymbol{\varphi}_{i}||_2^2)
\end{equation}
that
\begin{equation}
\nonumber
E(\boldsymbol{\varphi}_{i+1})-E(\boldsymbol{\varphi}_{i}) \approx - {\Delta s} ||\nabla E(\boldsymbol{\varphi}_{i})||_2^2 +\frac{1}{2} {\Delta s}^2\nabla E(\boldsymbol{\varphi}_{i})^T\nabla^2 E(\boldsymbol{\varphi}_{i}) \nabla E(\boldsymbol{\varphi}_{i}).
\end{equation}

To ensure stability $E(\boldsymbol{\varphi}_{i+1})-E(\boldsymbol{\varphi}_{i})\leq 0$,
we have

\begin{eqnarray}
\nabla E(\boldsymbol{\varphi}_{i})^T\big(\mathbf{I}-\frac{\Delta s}{2}\nabla^2 E(\boldsymbol{\varphi}_{i})\big) \nabla E(\boldsymbol{\varphi}_{i})\geq 0.
\label{ineq_1}
\end{eqnarray}
When $\nabla^2 E(\boldsymbol{\varphi}_{i})\leq 0$,  Eq.(\ref{ineq_1}) always holds. Otherwise $\Delta s$ needs to satisfy
\begin{eqnarray}
\nonumber
\Delta s\leq{\Delta s}_{stable}:={2}/{\rho^+\big(\nabla^2 E(\boldsymbol{\varphi}_{i})\big)}
\end{eqnarray}
where $\rho^+\big(\nabla^2 E(\boldsymbol{\varphi}_{i})\big)$ is the maximum of the positive eigenvalues of $\nabla^2 E(\boldsymbol{\varphi}_{i})$. 

Next we derive $\Delta s_{opt}$ from Eq.(\ref{optimal}). First we rewrite the iteration Eq.(\ref{iter_rule}) as \begin{eqnarray}
\nonumber
\boldsymbol{\varphi}_{i+1}-\boldsymbol{\varphi}_* &=& \boldsymbol{\varphi}_{i}-\boldsymbol{\varphi}_* -\Delta s \big( \nabla E(\boldsymbol{\varphi}_{i}) - \nabla E(\boldsymbol{\varphi}_*) \big)\\\nonumber
&=&\bigg(\mathbf{I}-\Delta s \int_0^1\nabla^2 E\big(\boldsymbol{\varphi}_{i}+r(\boldsymbol{\varphi}_*- \boldsymbol{\varphi}_{i})\big)dr \bigg)(\boldsymbol{\varphi}_{i}-\boldsymbol{\varphi}_*)\\\nonumber
&=&\bigg(\mathbf{I}-\Delta s \int_0^1\nabla^2 E(\boldsymbol{\varphi}_{i})dr \bigg)(\boldsymbol{\varphi}_{i}-\boldsymbol{\varphi}_*)\\\nonumber
&&-\bigg(\Delta s \int_0^1\big(\nabla^2 E\big(\boldsymbol{\varphi}_{i}+r(\boldsymbol{\varphi}_*- \boldsymbol{\varphi}_{i})\big)-\nabla^2 E(\boldsymbol{\varphi}_{i})\big)dr \bigg)(\boldsymbol{\varphi}_{i}-\boldsymbol{\varphi}_*).
\end{eqnarray}
where $\boldsymbol{\varphi}_*$ is a critical point such that $\nabla E(\boldsymbol{\varphi}_{*})=0$. So far we have only assumed that $\nabla E(\boldsymbol{\varphi}_{*})=0$. Following this, we have 
 \begin{eqnarray}
\nonumber
||\boldsymbol{\varphi}_{i+1}-\boldsymbol{\varphi}_*||_2&\leq&||\mathbf{I}-\Delta s \int_0^1\nabla^2 E(\boldsymbol{\varphi}_{i})dr ||_2||\boldsymbol{\varphi}_{i}-\boldsymbol{\varphi}_*||_2\\\nonumber
&&+||\Delta s \int_0^1\big(\nabla^2 E\big(\boldsymbol{\varphi}_{i}+r(\boldsymbol{\varphi}_*- \boldsymbol{\varphi}_{i})\big)-\nabla^2 E(\boldsymbol{\varphi}_{i})\big)dr ||_2||\boldsymbol{\varphi}_{i}-\boldsymbol{\varphi}_*||_2.
\end{eqnarray}

Further assuming  $\nabla^2 E$ is Lipschitz continuous with constant $L$, we have 
 \begin{eqnarray}
\nonumber
||\boldsymbol{\varphi}_{i+1}-\boldsymbol{\varphi}_*||_2&\leq&||\mathbf{I}-\Delta s \int_0^1\nabla^2 E(\boldsymbol{\varphi}_{i})dr ||_2||\boldsymbol{\varphi}_{i}-\boldsymbol{\varphi}_*||_2\\\nonumber
&&+\frac{\Delta sL}{2}  ||\boldsymbol{\varphi}_{i}-\boldsymbol{\varphi}_*||^2_2.
\end{eqnarray}
Thus, the rate of convergence $\lambda :=||\boldsymbol{\varphi}_{i+1}-\boldsymbol{\varphi}_* ||_2/||\boldsymbol{\varphi}_{i}-\boldsymbol{\varphi}_*||_2$ is bounded above by 
\begin{eqnarray}
\lambda&\leq& \rho\big(\mathbf{I}-\Delta s \nabla^2 E\big(\boldsymbol{\varphi}_{i}\big) \big) + \frac{\Delta sL}{2}  ||\boldsymbol{\varphi}_{i}-\boldsymbol{\varphi}_*||_2
\label{lip}
\end{eqnarray}
where $\rho(\Box)$ is the spectral radius of $\Box$.

Finally we assume $||\boldsymbol{\varphi}_i-\boldsymbol{\varphi}_*||_2\ll1$, the second term in Eq.(\ref{lip}) can be dropped:
\begin{eqnarray}
\nonumber\lambda\leq \lambda_{bound}(\Delta s)&=:&\rho\big(\mathbf{I}-\Delta s \nabla^2 E\big(\boldsymbol{\varphi}_{i}\big) \big)\\\nonumber
&=&\max{\bigg(|1-\Delta s\lambda_{min}\big(\nabla^2 E(\boldsymbol{\varphi}_{i})\big)|,|1-\Delta s\lambda_{max}\big(\nabla^2 E(\boldsymbol{\varphi}_{i})\big)|\bigg)},\end{eqnarray}
and our aim is to find $\Delta s_{opt} = \arg \min_{\Delta s} \lambda_{bound}(\Delta s)\geq 0$.

Notice if $\nabla^2 E(\boldsymbol{\varphi}_{i})$ is not positive semi-definite ($\lambda_{min}\big(\nabla^2 E(\boldsymbol{\varphi}_{i})\big)<0$), we have $ \lambda_{bound}(\Delta s)>1$ for any $\Delta s>0$. Thus we consider $\nabla^2E(\boldsymbol{\varphi}_{i})\geq0$   and find  \begin{equation}\Delta s_{opt} = {2}/\bigg(\lambda_{max}\big(\nabla^2 E(\boldsymbol{\varphi}_{i})\big)+ \lambda_{min}\big(\nabla^2 E(\boldsymbol{\varphi}_{i})\big)\bigg)\nonumber\end{equation} 
as the minimizer of $\lambda_{bound}(\Delta s)$ (Eq.(\ref{optimal}) in the main text), with the corresponding optimal rate \begin{equation}\lambda_{opt} = \frac{\lambda_{max}\big(\nabla^2 E(\boldsymbol{\varphi}_{i})\big)- \lambda_{min}\big(\nabla^2 E(\boldsymbol{\varphi}_{i})\big)}{\lambda_{max}\big(\nabla^2 E(\boldsymbol{\varphi}_{i})\big)+ \lambda_{min}\big(\nabla^2 E(\boldsymbol{\varphi}_{i})\big)}\nonumber.\end{equation}

Notice ${\Delta s}_{opt}$ is smaller than ${\Delta s}_{stable}$ from Eq. (\ref{stable}) when  $\nabla^2 E(\boldsymbol{\varphi}_{i})>0$. Thus we choose $\Delta {s} = \min{(\Delta {s}_{opt}, \Delta s_{stable})}$.

\section*{Acknowledgment}
This work is partially supported by the National Science Foundation under award number NSF-DMS-2012330 and  CAREER-2144372 and by the National Institute of General Medical Sciences of the National Institutes of Health under award number R01GM157590.
\bibliographystyle{elsarticle-harv}
\bibliography{elsarticle-template-harv_hidden}

\begin{thebibliography}{39}
\expandafter\ifx\csname natexlab\endcsname\relax\def\natexlab#1{#1}\fi
\providecommand{\url}[1]{\texttt{#1}}
\providecommand{\href}[2]{#2}
\providecommand{\path}[1]{#1}
\providecommand{\DOIprefix}{doi:}
\providecommand{\ArXivprefix}{arXiv:}
\providecommand{\URLprefix}{URL: }
\providecommand{\Pubmedprefix}{pmid:}
\providecommand{\doi}[1]{\href{http://dx.doi.org/#1}{\path{#1}}}
\providecommand{\Pubmed}[1]{\href{pmid:#1}{\path{#1}}}
\providecommand{\bibinfo}[2]{#2}
\ifx\xfnm\relax \def\xfnm[#1]{\unskip,\space#1}\fi
\bibitem[{Altorki et~al.(2019)Altorki, Markowitz, Gao, Port, Saxena, Stiles,
  McGraw and Mittal}]{altorki2019lung}
\bibinfo{author}{Altorki, N.K.}, \bibinfo{author}{Markowitz, G.J.},
  \bibinfo{author}{Gao, D.}, \bibinfo{author}{Port, J.L.},
  \bibinfo{author}{Saxena, A.}, \bibinfo{author}{Stiles, B.},
  \bibinfo{author}{McGraw, T.}, \bibinfo{author}{Mittal, V.},
  \bibinfo{year}{2019}.
\newblock \bibinfo{title}{The lung microenvironment: an important regulator of
  tumour growth and metastasis}.
\newblock \bibinfo{journal}{Nature Reviews Cancer} \bibinfo{volume}{19},
  \bibinfo{pages}{9--31}.
\bibitem[{Ambrosi et~al.(2019)Ambrosi, Ben~Amar, Cyron, DeSimone, Goriely,
  Humphrey and Kuhl}]{Ambrosi-2019-growth-review}
\bibinfo{author}{Ambrosi, D.}, \bibinfo{author}{Ben~Amar, M.},
  \bibinfo{author}{Cyron, C.J.}, \bibinfo{author}{DeSimone, A.},
  \bibinfo{author}{Goriely, A.}, \bibinfo{author}{Humphrey, J.D.},
  \bibinfo{author}{Kuhl, E.}, \bibinfo{year}{2019}.
\newblock \bibinfo{title}{{Growth and remodelling of living tissues:
  Perspectives, challenges and opportunities}}.
\newblock \bibinfo{journal}{Journal of the Royal Society Interface}
  \bibinfo{volume}{16}.
\newblock \DOIprefix\doi{10.1098/rsif.2019.0233}.
\bibitem[{Audoly and Pomeau(2000)}]{audoly2000elasticity}
\bibinfo{author}{Audoly, B.}, \bibinfo{author}{Pomeau, Y.},
  \bibinfo{year}{2000}.
\newblock \bibinfo{title}{Elasticity and geometry}, in:
  \bibinfo{booktitle}{Peyresq Lectures on Nonlinear Phenomena}.
  \bibinfo{publisher}{World Scientific}, pp. \bibinfo{pages}{1--35}.
\bibitem[{Balbi et~al.(2015)Balbi, Kuhl and Ciarletta}]{balbi2015morphoelastic}
\bibinfo{author}{Balbi, V.}, \bibinfo{author}{Kuhl, E.},
  \bibinfo{author}{Ciarletta, P.}, \bibinfo{year}{2015}.
\newblock \bibinfo{title}{Morphoelastic control of gastro-intestinal
  organogenesis: Theoretical predictions and numerical insights}.
\newblock \bibinfo{journal}{Journal of the Mechanics and Physics of Solids}
  \bibinfo{volume}{78}, \bibinfo{pages}{493--510}.
\bibitem[{Ben~Amar(2023)}]{amar2023creases}
\bibinfo{author}{Ben~Amar, M.}, \bibinfo{year}{2023}.
\newblock \bibinfo{title}{Creases and cusps in growing soft matter}.
\newblock \bibinfo{journal}{arXiv preprint arXiv:2309.11412} .
\bibitem[{Ben~Amar and Goriely(2005)}]{amar2005growth}
\bibinfo{author}{Ben~Amar, M.}, \bibinfo{author}{Goriely, A.},
  \bibinfo{year}{2005}.
\newblock \bibinfo{title}{Growth and instability in elastic tissues}.
\newblock \bibinfo{journal}{Journal of the Mechanics and Physics of Solids}
  \bibinfo{volume}{53}, \bibinfo{pages}{2284--2319}.
\bibitem[{Budday et~al.(2015)Budday, Steinmann, Goriely and
  Kuhl}]{budday2015size}
\bibinfo{author}{Budday, S.}, \bibinfo{author}{Steinmann, P.},
  \bibinfo{author}{Goriely, A.}, \bibinfo{author}{Kuhl, E.},
  \bibinfo{year}{2015}.
\newblock \bibinfo{title}{Size and curvature regulate pattern selection in the
  mammalian brain}.
\newblock \bibinfo{journal}{Extreme Mechanics Letters} \bibinfo{volume}{4},
  \bibinfo{pages}{193--198}.
\bibitem[{Ciarletta and Ben~Amar(2012)}]{ciarletta2012growth}
\bibinfo{author}{Ciarletta, P.}, \bibinfo{author}{Ben~Amar, M.},
  \bibinfo{year}{2012}.
\newblock \bibinfo{title}{Growth instabilities and folding in tubular organs: a
  variational method in non-linear elasticity}.
\newblock \bibinfo{journal}{International Journal of Non-Linear Mechanics}
  \bibinfo{volume}{47}, \bibinfo{pages}{248--257}.
\bibitem[{Folkman(2002)}]{folkman2002role}
\bibinfo{author}{Folkman, J.}, \bibinfo{year}{2002}.
\newblock \bibinfo{title}{Role of angiogenesis in tumor growth and metastasis},
  in: \bibinfo{booktitle}{Seminars in oncology},
  \bibinfo{organization}{Elsevier}. pp. \bibinfo{pages}{15--18}.
\bibitem[{Garcia et~al.(2017)Garcia, Okamoto, Bayly and
  Taber}]{garcia2017contraction}
\bibinfo{author}{Garcia, K.E.}, \bibinfo{author}{Okamoto, R.J.},
  \bibinfo{author}{Bayly, P.V.}, \bibinfo{author}{Taber, L.A.},
  \bibinfo{year}{2017}.
\newblock \bibinfo{title}{Contraction and stress-dependent growth shape the
  forebrain of the early chicken embryo}.
\newblock \bibinfo{journal}{Journal of the mechanical behavior of biomedical
  materials} \bibinfo{volume}{65}, \bibinfo{pages}{383--397}.
\bibitem[{Goriely(2017)}]{goriely2017mathematics}
\bibinfo{author}{Goriely, A.}, \bibinfo{year}{2017}.
\newblock \bibinfo{title}{The mathematics and mechanics of biological growth}.
\newblock \bibinfo{publisher}{Springer}.
\bibitem[{Goriely and Ben~Amar(2005)}]{goriely2005differential}
\bibinfo{author}{Goriely, A.}, \bibinfo{author}{Ben~Amar, M.},
  \bibinfo{year}{2005}.
\newblock \bibinfo{title}{Differential growth and instability in elastic
  shells}.
\newblock \bibinfo{journal}{Physical review letters} \bibinfo{volume}{94},
  \bibinfo{pages}{198103}.
\bibitem[{Groh(2022)}]{groh2022morphoelastic}
\bibinfo{author}{Groh, R.M.}, \bibinfo{year}{2022}.
\newblock \bibinfo{title}{A morphoelastic stability framework for post-critical
  pattern formation in growing thin biomaterials}.
\newblock \bibinfo{journal}{Computer Methods in Applied Mechanics and
  Engineering} \bibinfo{volume}{394}, \bibinfo{pages}{114839}.
\bibitem[{Hallatschek et~al.(2023)Hallatschek, Datta, Drescher, Dunkel, Elgeti,
  Waclaw and Wingreen}]{hallatschek2023proliferating}
\bibinfo{author}{Hallatschek, O.}, \bibinfo{author}{Datta, S.S.},
  \bibinfo{author}{Drescher, K.}, \bibinfo{author}{Dunkel, J.},
  \bibinfo{author}{Elgeti, J.}, \bibinfo{author}{Waclaw, B.},
  \bibinfo{author}{Wingreen, N.S.}, \bibinfo{year}{2023}.
\newblock \bibinfo{title}{Proliferating active matter}.
\newblock \bibinfo{journal}{Nature Reviews Physics} \bibinfo{volume}{5},
  \bibinfo{pages}{407--419}.
\bibitem[{Hofer and Lutolf(2021)}]{hofer2021engineering}
\bibinfo{author}{Hofer, M.}, \bibinfo{author}{Lutolf, M.P.},
  \bibinfo{year}{2021}.
\newblock \bibinfo{title}{Engineering organoids}.
\newblock \bibinfo{journal}{Nature Reviews Materials} \bibinfo{volume}{6},
  \bibinfo{pages}{402--420}.
\bibitem[{Holzapfel(2002)}]{holzapfel2002nonlinear}
\bibinfo{author}{Holzapfel, G.A.}, \bibinfo{year}{2002}.
\newblock \bibinfo{title}{Nonlinear solid mechanics: a continuum approach for
  engineering science}.
\bibitem[{Ingber(2006)}]{ingber2006mechanical}
\bibinfo{author}{Ingber, D.E.}, \bibinfo{year}{2006}.
\newblock \bibinfo{title}{Mechanical control of tissue morphogenesis during
  embryological development}.
\newblock \bibinfo{journal}{The International journal of developmental biology}
  .
\bibitem[{Jin et~al.(2011)Jin, Cai and Suo}]{jin2011creases}
\bibinfo{author}{Jin, L.}, \bibinfo{author}{Cai, S.}, \bibinfo{author}{Suo,
  Z.}, \bibinfo{year}{2011}.
\newblock \bibinfo{title}{Creases in soft tissues generated by growth}.
\newblock \bibinfo{journal}{Europhysics Letters} \bibinfo{volume}{95},
  \bibinfo{pages}{64002}.
\bibitem[{Kr{\"o}ner(1959)}]{kroner1959allgemeine}
\bibinfo{author}{Kr{\"o}ner, E.}, \bibinfo{year}{1959}.
\newblock \bibinfo{title}{Allgemeine kontinuumstheorie der versetzungen und
  eigenspannungen}.
\newblock \bibinfo{journal}{Archive for Rational Mechanics and Analysis}
  \bibinfo{volume}{4}, \bibinfo{pages}{273--334}.
\bibitem[{Lancaster and Knoblich(2014)}]{lancaster2014organogenesis}
\bibinfo{author}{Lancaster, M.A.}, \bibinfo{author}{Knoblich, J.A.},
  \bibinfo{year}{2014}.
\newblock \bibinfo{title}{Organogenesis in a dish: modeling development and
  disease using organoid technologies}.
\newblock \bibinfo{journal}{Science} \bibinfo{volume}{345},
  \bibinfo{pages}{1247125}.
\bibitem[{Lecuit and Pilot(2003)}]{lecuit2003developmental}
\bibinfo{author}{Lecuit, T.}, \bibinfo{author}{Pilot, F.},
  \bibinfo{year}{2003}.
\newblock \bibinfo{title}{Developmental control of cell morphogenesis: a focus
  on membrane growth}.
\newblock \bibinfo{journal}{Nature cell biology} \bibinfo{volume}{5},
  \bibinfo{pages}{103--108}.
\bibitem[{Lee(1969)}]{lee1969elastic}
\bibinfo{author}{Lee, E.H.}, \bibinfo{year}{1969}.
\newblock \bibinfo{title}{Elastic-plastic deformation at finite strains}.
\newblock \bibinfo{journal}{Trans.ASMEJ.Appl.Mech.} \bibinfo{volume}{54},
  \bibinfo{pages}{1--6}.
\bibitem[{Li et~al.(2011a)Li, Cao, Feng and Gao}]{li2011surface}
\bibinfo{author}{Li, B.}, \bibinfo{author}{Cao, Y.P.}, \bibinfo{author}{Feng,
  X.Q.}, \bibinfo{author}{Gao, H.}, \bibinfo{year}{2011}a.
\newblock \bibinfo{title}{Surface wrinkling of mucosa induced by volumetric
  growth: theory, simulation and experiment}.
\newblock \bibinfo{journal}{Journal of the Mechanics and Physics of Solids}
  \bibinfo{volume}{59}, \bibinfo{pages}{758--774}.
\bibitem[{Li et~al.(2011b)Li, Jia, Cao, Feng and Gao}]{li2011surface2}
\bibinfo{author}{Li, B.}, \bibinfo{author}{Jia, F.}, \bibinfo{author}{Cao,
  Y.P.}, \bibinfo{author}{Feng, X.Q.}, \bibinfo{author}{Gao, H.},
  \bibinfo{year}{2011}b.
\newblock \bibinfo{title}{Surface wrinkling patterns on a core-shell soft
  sphere}.
\newblock \bibinfo{journal}{Physical review letters} \bibinfo{volume}{106},
  \bibinfo{pages}{234301}.
\bibitem[{Loewenstein and Penn(1967)}]{loewenstein1967intercellular}
\bibinfo{author}{Loewenstein, W.R.}, \bibinfo{author}{Penn, R.D.},
  \bibinfo{year}{1967}.
\newblock \bibinfo{title}{Intercellular communication and tissue growth: Ii.
  tissue regeneration}.
\newblock \bibinfo{journal}{The Journal of Cell Biology} \bibinfo{volume}{33},
  \bibinfo{pages}{235--242}.
\bibitem[{Moulton and Goriely(2011)}]{moulton2011circumferential}
\bibinfo{author}{Moulton, D.}, \bibinfo{author}{Goriely, A.},
  \bibinfo{year}{2011}.
\newblock \bibinfo{title}{Circumferential buckling instability of a growing
  cylindrical tube}.
\newblock \bibinfo{journal}{Journal of the Mechanics and Physics of Solids}
  \bibinfo{volume}{59}, \bibinfo{pages}{525--537}.
\bibitem[{Neufeld et~al.(1998)Neufeld, De~La~Cruz, Johnston and
  Edgar}]{neufeld1998coordination}
\bibinfo{author}{Neufeld, T.P.}, \bibinfo{author}{De~La~Cruz, A.F.A.},
  \bibinfo{author}{Johnston, L.A.}, \bibinfo{author}{Edgar, B.A.},
  \bibinfo{year}{1998}.
\newblock \bibinfo{title}{Coordination of growth and cell division in the
  drosophila wing}.
\newblock \bibinfo{journal}{Cell} \bibinfo{volume}{93},
  \bibinfo{pages}{1183--1193}.
\bibitem[{Nia et~al.(2016)Nia, Liu, Seano, Datta, Jones, Rahbari, Incio,
  Chauhan, Jung, Martin et~al.}]{nia2016solid}
\bibinfo{author}{Nia, H.T.}, \bibinfo{author}{Liu, H.}, \bibinfo{author}{Seano,
  G.}, \bibinfo{author}{Datta, M.}, \bibinfo{author}{Jones, D.},
  \bibinfo{author}{Rahbari, N.}, \bibinfo{author}{Incio, J.},
  \bibinfo{author}{Chauhan, V.P.}, \bibinfo{author}{Jung, K.},
  \bibinfo{author}{Martin, J.D.}, et~al., \bibinfo{year}{2016}.
\newblock \bibinfo{title}{Solid stress and elastic energy as measures of tumour
  mechanopathology}.
\newblock \bibinfo{journal}{Nature biomedical engineering} \bibinfo{volume}{1},
  \bibinfo{pages}{0004}.
\bibitem[{Oltean and Taber(2018)}]{oltean2018apoptosis}
\bibinfo{author}{Oltean, A.}, \bibinfo{author}{Taber, L.A.},
  \bibinfo{year}{2018}.
\newblock \bibinfo{title}{Apoptosis generates mechanical forces that close the
  lens vesicle in the chick embryo}.
\newblock \bibinfo{journal}{Physical biology} \bibinfo{volume}{15},
  \bibinfo{pages}{025001}.
\bibitem[{Papastavrou et~al.(2013)Papastavrou, Steinmann and
  Kuhl}]{papastavrou2013mechanics}
\bibinfo{author}{Papastavrou, A.}, \bibinfo{author}{Steinmann, P.},
  \bibinfo{author}{Kuhl, E.}, \bibinfo{year}{2013}.
\newblock \bibinfo{title}{On the mechanics of continua with boundary energies
  and growing surfaces}.
\newblock \bibinfo{journal}{Journal of the Mechanics and Physics of Solids}
  \bibinfo{volume}{61}, \bibinfo{pages}{1446--1463}.
\bibitem[{Podio-Guidugli and Caffarelli(1991)}]{podio1991surface}
\bibinfo{author}{Podio-Guidugli, P.}, \bibinfo{author}{Caffarelli, G.V.},
  \bibinfo{year}{1991}.
\newblock \bibinfo{title}{Surface interaction potentials in elasticity}, in:
  \bibinfo{booktitle}{Mechanics and Thermodynamics of Continua: A Collection of
  Papers Dedicated to BD Coleman on His Sixtieth Birthday}.
  \bibinfo{publisher}{Springer}, pp. \bibinfo{pages}{345--385}.
\bibitem[{Reddy(2017)}]{reddy2017energy}
\bibinfo{author}{Reddy, J.N.}, \bibinfo{year}{2017}.
\newblock \bibinfo{title}{Energy principles and variational methods in applied
  mechanics}.
\newblock \bibinfo{publisher}{John Wiley \& Sons}.
\bibitem[{Rodriguez et~al.(1994)Rodriguez, Hoger and
  McCulloch}]{Rodriguez-1994-decomposition}
\bibinfo{author}{Rodriguez, E.K.}, \bibinfo{author}{Hoger, A.},
  \bibinfo{author}{McCulloch, A.D.}, \bibinfo{year}{1994}.
\newblock \bibinfo{title}{{Stress-dependent finite growth in soft elastic
  tissues}}.
\newblock \bibinfo{journal}{Journal of Biomechanics} \bibinfo{volume}{27},
  \bibinfo{pages}{455--467}.
\newblock \DOIprefix\doi{10.1016/0021-9290(94)90021-3}.
\bibitem[{Sun and Irvine(2014)}]{sun2014control}
\bibinfo{author}{Sun, G.}, \bibinfo{author}{Irvine, K.D.},
  \bibinfo{year}{2014}.
\newblock \bibinfo{title}{Control of growth during regeneration}.
\newblock \bibinfo{journal}{Current topics in developmental biology}
  \bibinfo{volume}{108}, \bibinfo{pages}{95--120}.
\bibitem[{Taber(1995)}]{taber1995biomechanics}
\bibinfo{author}{Taber, L.}, \bibinfo{year}{1995}.
\newblock \bibinfo{title}{Biomechanics of growth, remodeling, and
  morphogenesis}.
\newblock \bibinfo{journal}{Journal of Biomechanical Engineering}
  \bibinfo{volume}{117}, \bibinfo{pages}{429–439}.
\newblock \DOIprefix\doi{10.1115/1.2794191}.
\bibitem[{Tallinen et~al.(2013)Tallinen, Biggins and
  Mahadevan}]{tallinen2013surface}
\bibinfo{author}{Tallinen, T.}, \bibinfo{author}{Biggins, J.S.},
  \bibinfo{author}{Mahadevan, L.}, \bibinfo{year}{2013}.
\newblock \bibinfo{title}{Surface sulci in squeezed soft solids}.
\newblock \bibinfo{journal}{Physical review letters} \bibinfo{volume}{110},
  \bibinfo{pages}{024302}.
\bibitem[{Tallinen et~al.(2016)Tallinen, Chung, Rousseau, Girard, Lef{\`e}vre
  and Mahadevan}]{tallinen2016growth}
\bibinfo{author}{Tallinen, T.}, \bibinfo{author}{Chung, J.Y.},
  \bibinfo{author}{Rousseau, F.}, \bibinfo{author}{Girard, N.},
  \bibinfo{author}{Lef{\`e}vre, J.}, \bibinfo{author}{Mahadevan, L.},
  \bibinfo{year}{2016}.
\newblock \bibinfo{title}{On the growth and form of cortical convolutions}.
\newblock \bibinfo{journal}{Nature Physics} \bibinfo{volume}{12},
  \bibinfo{pages}{588--593}.
\bibitem[{Wu and {Ben Amar}(2015)}]{Wu-2015-wound}
\bibinfo{author}{Wu, M.}, \bibinfo{author}{{Ben Amar}, M.},
  \bibinfo{year}{2015}.
\newblock \bibinfo{title}{{Growth and remodelling for profound circular wounds
  in skin}}.
\newblock \bibinfo{journal}{Biomechanics and Modeling in Mechanobiology}
  \bibinfo{volume}{14}, \bibinfo{pages}{357--370}.
\newblock \URLprefix \url{http://link.springer.com/10.1007/s10237-014-0609-1},
  \DOIprefix\doi{10.1007/s10237-014-0609-1}.
\bibitem[{Wu and Ben~Amar(2015)}]{wu2015modelling}
\bibinfo{author}{Wu, M.}, \bibinfo{author}{Ben~Amar, M.}, \bibinfo{year}{2015}.
\newblock \bibinfo{title}{Modelling fibers in growing disks of soft tissues}.
\newblock \bibinfo{journal}{Mathematics and Mechanics of Solids}
  \bibinfo{volume}{20}, \bibinfo{pages}{663--679}.

\end{thebibliography}

\end{document}